\documentclass[11pt, dina4,  amssymb]{article}
\usepackage{amssymb, amsmath,amsthm}
\usepackage[titletoc, title]{appendix}
\usepackage{graphicx}
\usepackage{diagbox}
\usepackage{cite}
\usepackage{multirow}
\usepackage{epstopdf}
\usepackage{color}
\usepackage[margin=1in]{geometry}
\usepackage[normalem]{ulem}

\newcommand{\p}{\partial}
\newcommand{\og}{\omega}
\newcommand{\Og}{\Omega}
\newcommand{\fl}[2]{\frac{#1}{#2}}

\newcommand{\nn}{\nonumber}
\newcommand{\ap}{\alpha}

\newcommand{\veps}{\varepsilon}

\newcommand{\Dt}{\Delta}

\newcommand{\be}{\begin{equation}}
\newcommand{\ee}{\end{equation}}
\newcommand{\ba}{\begin{array}}
\newcommand{\ea}{\end{array}}
\newcommand{\bea}{\begin{eqnarray}}
\newcommand{\eea}{\end{eqnarray}}
\newcommand{\beas}{\begin{eqnarray*}}
\newcommand{\eeas}{\end{eqnarray*}}
\newtheorem{remark}{Remark}[section]
\newtheorem{lemma}{Lemma}[section]

\newcommand{\bx}{{\bf x} }
\newcommand{\by}{{\bf y} }

\newcommand{\bb}{\vskip 10pt}
\definecolor{ForestGreen}{rgb}{0.0, 0.5, 0.0}

\title{A unified meshfree pseudospectral method for solving both classical and fractional PDEs}
\author{John Burkardt\thanks{Department of Mathematics, University of Pittsburgh, Pittsburgh, PA 15260 (Email: jvburkardt@gmail.com)},  \ \
Yixuan Wu\thanks{Department of Mathematics and Statistics, Missouri University of Science and Technology, Rolla, MO 65409 (Email:  yxw7c@mst.edu)}, \ \
Yanzhi Zhang\thanks{Department of Mathematics and Statistics, Missouri University of Science and Technology, Rolla, MO 65409 (Email:  zhangyanz@mst.edu)}
}
\begin{document}
\date{}
\maketitle

\begin{abstract} 
In this paper, we propose a meshfree  method based on the Gaussian radial basis function (RBF) to solve both classical and fractional PDEs. 
The proposed method takes advantage of the analytical  Laplacian of Gaussian functions so as to accommodate the discretization of the classical and fractional Laplacian in a single framework and avoid the large computational cost for numerical evaluation of the fractional derivatives. 
These important merits distinguish it from other numerical methods for fractional PDEs. 
Moreover, our method is simple and easy to handle complex geometry and local refinement, and its computer program implementation remains the same for any dimension $d \ge 1$.  
Extensive numerical experiments are provided to study the performance of our method in both approximating the Dirichlet Laplace operators and solving PDE problems. 
Compared to the recently proposed Wendland RBF method, our method  exactly incorporates the Dirichlet boundary conditions into the scheme and is free of the Gibbs phenomenon as observed in the literature. 
Our studies suggest that to obtain good accuracy the shape parameter cannot be too small or too big, and the optimal shape parameter might depend on the RBF center points and the solution properties. 
\end{abstract}


\section{Introduction}
\setcounter{equation}{0}
\label{section1}

In the recent decade, fractional partial differential equations (PDEs) have found widespread applications in many fields, including turbulence \cite{Del2004, Epps0018, Egolf}, geophysics \cite{Baeumer2010, Zhang2012}, biomedicine and biology \cite{Magin2009, Javanainen2013}, and quantum mechanics \cite{Laskin2000, Duo2015}. 
In traditional (integer-order) PDEs, diffusion describes the transport process due to Brownian motion and  is modeled by the classical Laplacian $\Dt$. 
In contrast, diffusion in fractional PDEs is described by the fractional Laplacian $(-\Dt)^\fl{\ap}{2}$ which underlines the L\'evy transport. 
It was recently found that  the Brownian and L\'evy transports might coexist in many complex (e.g., biological and chemical) systems  \cite{Baeumer2010, Zhang2012, Javanainen2013}. 
Thus mathematical models including both classical and fractional Laplacians could be more proper to describe such a phenomenon. 
On the other hand, the classical and fractional Laplacians, one local and the other nonlocal,  possess  distinct properties. 
Consequently, the analytical and numerical frameworks for studying  these two operators are significantly different. 
For instance,  numerical discretizations (e.g., finite element methods) for the classical and fractional Laplacians are usually incompatible, and separate implementation efforts are required to study problems of these two operators.  
In this work, we propose a new meshfree pseudospectral method with the intrinsic merit of solving  both  classical and fractional PDEs in a unified scheme.

The classical and fractional Laplacians can be defined via the  parametric pseudo-differential operator with symbol $|\xi|^\ap$ \cite{Landkof,Samko}:
\begin{eqnarray}
\label{pseudo}
(-\Delta)^{\fl{\alpha}{2}}u({\bx}) = \mathcal{F}^{-1}\big[|\xi|^\alpha \mathcal{F}[u]\big], \qquad \mbox{for} \ \  \ap \ge 0, 
\end{eqnarray}
where $\mathcal{F}$ is the Fourier transform with associated inverse transform $\mathcal{F}^{-1}$. 
The definition (\ref{pseudo}) covers a wide class of operators for different values of the parameter $\ap$.
In this work, we are interested in the exponent $\ap \in (0, 2]$.
For $\ap = 2$,   the formulation \eqref{pseudo} gives the spectral representation of the {\it classical Laplacian} $-\Dt$,  while it is referred to as the {\it fractional Laplacian} if $\ap < 2$.  
Probabilistically, the fractional Laplacian  represents the infinitesimal generator of a symmetric $\ap$-stable L\'evy process. 
It can be also defined  in a hypersingular integral form (also known as the integral fractional Laplacian) \cite{Landkof,Samko,Kwasnicki2017}: 
\bea\label{integralFL}
(-\Delta)^{\fl{\alpha}{2}}u(\bx) = C_{d,\alpha}\,{\rm P. V.}\int_{{\mathbb R}^d}\fl{u(\bx) - u({\bf y})}{|\bx -{\bf y}|^{d+\ap}}d{\bf y}, \qquad \mbox{for} \ \  \ap \in (0, 2),
\eea
for  $d  = 1, 2$, or $3$, where ${\rm P. V.}$ stands for the principal value integral, and $|\bx - \by|$ denotes the Euclidean distance between points $\bx$ and $\by$. 
The normalization constant is given by
\beas
C_{d,\alpha}=\fl{2^{\ap-1} \ap\,\Gamma(({\ap+d})/{2})}{\sqrt{\pi^{d}}\,\Gamma(1 -{\ap}/{2})}
\eeas
with $\Gamma(\cdot)$ being the Gamma function. 
Over the entire space ${\mathbb R}^d$,
the integral fractional Laplacian (\ref{integralFL}) is equivalent to the pseudo-differential operator (\ref{pseudo}) with $\ap \in (0, 2)$ \cite{Kwasnicki2017,Samko,Landkof}.

The formulation in (\ref{pseudo}) provides a uniform definition of the classical and fractional Laplacians via the parametric symbol $|\xi|^\ap$. 
It suggests that if the entire space ${\mathbb R}^d$ or periodic bounded domains are considered, one can study these two operators together. 
For example, the Fourier pseudospectral methods based on (\ref{pseudo}) were introduced in  \cite{Duo2016, Kirkpatrick2016} to solve the classical and fractional Schr\"odinger equations on a periodic domain. 
However, if a non-periodic bounded domain is considered,  the pseudo-differential form of the classical and fractional Laplacians loses its advantages and has challenges to incorporate general boundary conditions.  
Thus different representations of the classical Laplacian (i.e. $\Dt = \p_{xx} + \p_{yy} + \p_{zz}$) and fractional Laplacian (i.e. formulation in (\ref{integralFL})) are adopted, which clearly manifests the differences between these two operators -- one is a local derivative operator, and the other is a nonlocal integral operator. 
In practice, numerical methods (e.g., finite difference/element methods) for these two operators  on domains with non-periodic boundary conditions are separately developed and  incompatible.

Compared to the classical Laplacian, numerical methods for  the fractional Laplacian (\ref{integralFL}) still remain limited. 
In \cite{Duo2018, Duo-FDM2019, Duo-TFL2019},  second-order finite difference methods were proposed to discretize the integral fractional Laplacian (\ref{integralFL}) for $d \ge 1$, and fast algorithms via the fast Fourier transforms were introduced for their efficient simulations. 
Various finite element  methods based on different formulations of the fractional Laplacian were developed in \cite{Acosta2017, Acosta2017-2D, Bonito2019, Ainsworth2018, DElia2013} to solve fractional problems. 
Recently, spectral methods were proposed to solve fractional PDEs with the integral fractional Laplacian  in bounded and unbounded domains \cite{Zhang2019, Tang2020}. 
These spectral methods could achieve higher accuracy than finite difference/element methods, but they have limited usability on irregular domains. 
So far the existing numerical methods for the classical Laplacian and their computer implementations cannot be used to solve problems with the fractional Laplacian, due to the distinct features of these two operators.

On the other hand, meshfree methods based on radial basis functions (RBF) have been widely applied to solve classical PDEs \cite{Fornberg, Fasshauer}. 
Compared to the mesh-based methods, these methods have more flexibility of domain geometry and can achieve higher accuracy with less computational cost.  
The application of RBF-based methods to solve fractional PDEs and nonlocal problems is still very recent. 
In \cite{Bond2015,Lehoucq2016,Zhao2018}, the Galerkin methods using a localized basis of RBFs were proposed to solve nonlocal diffusion problems. 
In \cite{Piret2013}, RBF-QR methods were proposed to solve the Riemann--Liouville spatial fractional diffusion problems. 
A Kansa RBF method was proposed in \cite{Pang2015} to solve the fractional advection-dispersion equations, where the spatial derivative was defined via  the fractional directional derivatives. 
Later, RBF collocation methods were introduced in \cite{Zafarghandi2019}  to solve similar advection-dispersion equations but with the Riesz spatial fractional derivatives. 
We remark that these RBF methods are for different fractional derivatives, and in this work we are interested in the fractional Laplacian. 
Recently, a Wendland RBF collocation method was proposed in \cite{Rosenfeld2019} to solve fractional problems with the  fractional Laplacian (\ref{pseudo}), { while a singular boundary method  based on a new definition of the fractional Laplacian was introduced in \cite{Chen2016}.}
Note that all the above RBF methods developed in the  fractional cases cannot be used to solve classical problems.

In this work, we propose a novel meshfree pseudospectral method based on the Gaussian RBFs, which has fundamental differences  from other RBF-based methods in \cite{Bond2015,Lehoucq2016,Zhao2018, Piret2013, Pang2015, Zafarghandi2019, Rosenfeld2019}. 
Inheriting the advantages of RBF methods, our method is simple and flexible of domain geometry, and its computer implementation remains the same for any dimension $d \ge 1$. 
Furthermore, it allows easy local refinements. 
Besides these advantages, our method has the distinct merits as below. \vspace{-2mm}
\begin{itemize}\itemsep -1pt
\item[(i)] It solves the classical and fractional PDEs in a unified scheme.  
To the best of our knowledge, this is the first numerical method that discretizes the classical and fractional Laplacians on non-periodic domains with a single scheme. 
This feature distinguishes our method from other existing methods which solve classical and fractional problems separately.
\item[(ii)] It takes great advantage of the Laplacian of the Gaussian RBFs (i.e. the confluent hypergeometric function) and avoids large computational costs in approximating the fractional derivative of RBFs, which is one major difference from those in \cite{Piret2013, Pang2015, Zafarghandi2019, Rosenfeld2019}. 
The fractional derivatives are usually defined in integral form  with a singular kernel. 
As pointed out in \cite{Pang2015}, it is challenging to balance the accuracy and efficiency in  approximating the fractional derivatives of RBFs with quadrature rules (e.g., Gauss--Jacobi quadrature rules are used in \cite{Pang2015}).
\item[(iii)] It exactly incorporates the boundary conditions into the scheme and  is free of the Runge phenomenon observed in \cite{Rosenfeld2019}. 
In contrast to the method  in \cite{Rosenfeld2019}, our method uses the exact boundary conditions and avoids evaluating the fractional Laplacian of RBFs with numerical quadrature rules. 
Moreover,  the method in \cite{Rosenfeld2019} requires a larger computational domain (much bigger than the actual physical domain), which significantly increases the computational costs especially for $d > 1$.
\end{itemize}
\vspace{-2mm}
The paper is organized as follows.  
In Section \ref{section2}, we first outline some important properties of the Gaussian RBFs and then introduce our numerical method on a bounded domain with Dirichlet boundary conditions. 
The performance of our method in approximating the Laplace operators is studied in  Section \ref{section3}. 
We then apply it solve the classical and fractional PDE problems  in  Section \ref{section4}. 
Finally,  some discussion and summary are made in Section \ref{section5}.

\section{RBF meshfree method}
\setcounter{equation}{0}
\label{section2}

Radial basis functions (RBFs) are well-known for their advantages in high-dimensional scattered data approximations and have recently broadened their applications in many areas, ranging from meteorology, statistics,  to machine learning. 
RBFs are usually real-valued scalar functions defined in the form of  $\varphi(|\bx|)$, where $|\bx|$ denotes the Euclidean norm of vector $\bx  \in  {\mathbb R}^d$. 
This radial form (e.g., $r = |\bx|$) makes their use for high-dimensional reconstruction problems very efficient and also allows invariance under orthogonal transforms. 
Common choices of RBFs usually fall into two main categories:  globally-supported functions (e.g., Gaussian RBFs), and compactly-supported functions (e.g.,  Wendland RBFs). 
For more discussion of RBFs, we refer the reader to \cite{Fornberg, Fasshauer} and references therein. 
In this work, we will use the Gaussian RBFs and start with some of their properties in Section \ref{section2-1}. 

\subsection{Gaussian radial basis functions}
\label{section2-1} 

Among all radial basis functions, the Gaussian RBF is a representative member of the class of infinitely differentiable functions with global support. 
It is defined as 
\beas\label{Gaussian}
\varphi(|\bx|) = \exp(-\varepsilon^2|\bx|^2), \qquad\mbox{for} \ \ \bx \in {\mathbb R}^d,
\eeas
where $\varepsilon \in {\mathbb R}$ denotes the shape parameter.  
The shape parameter $\varepsilon$ plays an important role in the approximation accuracy with Gaussian RBFs. 
It is usually chosen to be a constant, and recently spatial-dependent shape parameters were  studied in the literature \cite{Fornberg2017}.

When using RBF-based methods to solve fractional PDEs, the main challenge is to compute the fractional derivatives of RBFs \cite{Pang2015, Rosenfeld2019, Piret2013}. 
Their analytical solutions are usually unavailable, so numerical approximations are required to evaluate these fractional derivatives, e.g., the Gauss--Jacobi quadrature rules were used in \cite{Pang2015}. 
Note that the fractional derivatives are generally defined as an integral with singular kernel  over a large domain. 
Therefore, using quadrature rules to approximate the fractional derivatives of RBFs significantly increases the computational costs, especially in high-dimensional problems. 
Moreover,  it makes the implementation of RBF-based methods more complicated as special treatments are required around the singularity. 
These complications greatly deteriorate the performance of RBF-based methods in practice. 
In contrast, our method takes advantage of the properties of the Laplace operators and Gaussian RBFs so as to avoid numerical evaluations of the fractional derivatives with quadrature rules, which is one fundamental difference between our method and those in the literature \cite{Pang2015, Rosenfeld2019, Piret2013}.

To introduce our method, we will first present some important properties of the Laplace operators and their actions on Gaussian functions in the following lemmas. 
\begin{lemma}[{\bf  The Laplacian of Gaussian functions}]
\label{lemma1}
Let $u$ be a Gaussian function of the form  $u(\bx) = \exp(-|\bx|^2)$, for $\bx \in {\mathbb R}^d$. 
Then the Laplacian of $u$ is analytically given by \cite{Prudnikov, Dyda2017}:
\bea\label{Eq-lemma1}
(-\Dt)^{\fl{\ap}{2}} u(\bx) =  \fl{2^\ap\Gamma((d+\ap)/2)}{\Gamma(d/2)}\,_1F_1\Big(\fl{d+\ap}{2}; \, \fl{d}{2}; \, -|\bx|^2\Big), \qquad \mbox{for} \ \ \bx \in {\mathbb R}^d,\quad \ap \ge 0, 
\eea
where $_1F_1$ denotes the confluent hypergeometric function. 
\end{lemma}
Lemma \ref{lemma1} holds for any exponent $\ap \ge 0$. 
It  provides the foundation of developing unified schemes for the classical and fractional Laplacians. 
In the special case of $\ap = 2m$ with $m \in {\mathbb N}$, the result in (\ref{Eq-lemma1}) collapses to the classical integer-order derivatives $(-\Dt)^me^{-|\bx|^2}$. 
For instance,  using the properties of the confluent hypergeometric function $ _1F_1$, we obtain that (\ref{Eq-lemma1}) is equivalent to
\beas
-\Dt e^{-|\bx|^2} = e^{-|\bx|^2}(4|\bx|^2 - 2d), \qquad \mbox{for} \ \bx \in {\mathbb R}^d,
\eeas
when $\ap = 2$, i.e., the classical negative Laplacian of the Gaussian function.  

\begin{figure}[htb!]
\centerline{\includegraphics[height = 5.76cm, width = 7.86cm]{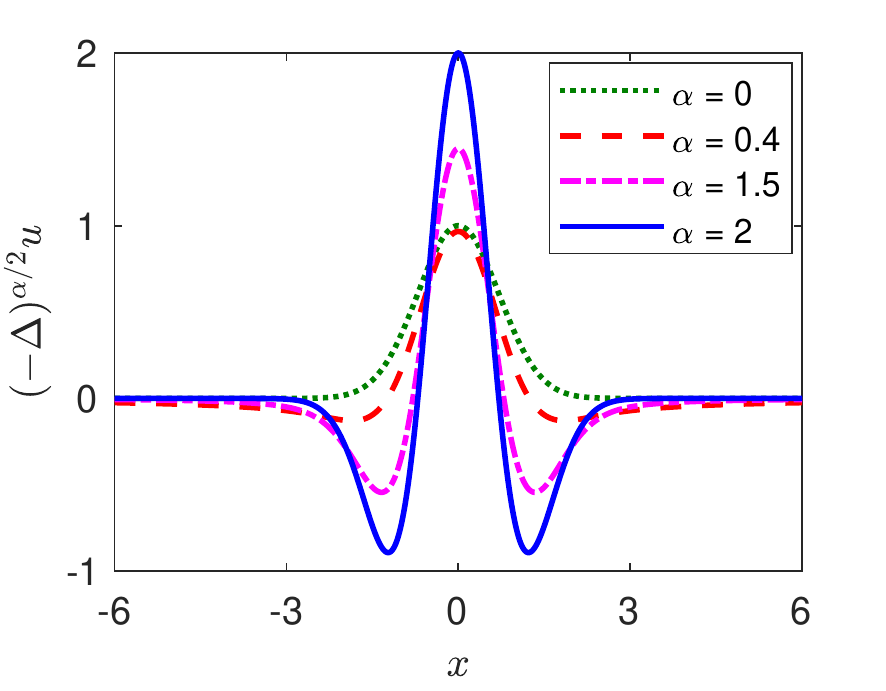}
\hspace{-3mm}
\includegraphics[height = 5.76cm, width = 7.86cm]{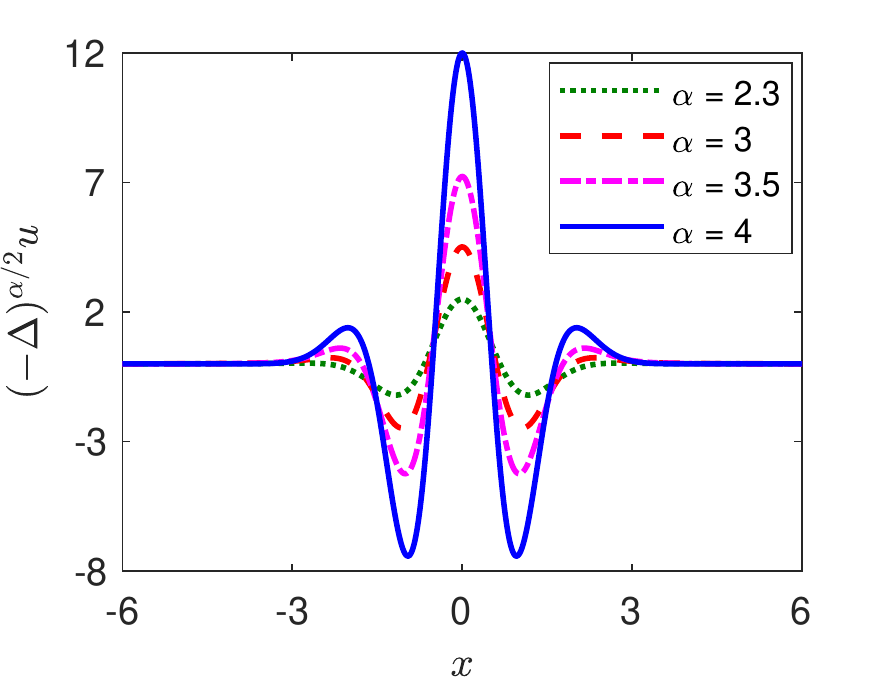}}
\caption{Illustration of function $(-\Dt)^\fl{\ap}{2}u$ with $u(x) = \exp(-x^2)$ for various $\ap \ge 0$.  }\label{Sec2-Fig1}
\end{figure}

Fig. \ref{Sec2-Fig1} illustrates the result $(-\Dt)^{\fl{\ap}{2}}u$ for various $\ap \ge 0$, where $u(x) = e^{-x^2}$ for $x\in{\mathbb R}$. 
It shows that the  function $(-\Dt)^\fl{\ap}{2}u$ is ``radially" symmetric with respect to $x = 0$, confirming that the Laplace operator is rotationally invariant. 
The rotational invariance is a crucial property in modeling isotropic anomalous diffusion in many applications \cite{Duo-FDM2019}. 
Moreover, the solution decays to zero as $|\bx| \to \infty$ -- the larger the exponent $\ap$, the faster the decay. 
In the extreme case of $\ap = 0$, there is  $(-\Dt)^{\fl{\ap}{2}}u(\bx) = u(\bx)$, that is, $(-\Dt)^{\fl{\ap}{2}}$ reduces to the identify operator $I$, consistent with the definition (\ref{pseudo}) as $\ap \to 0$. 
Fig. \ref{Sec2-Fig1} additionally shows that the function $(-\Dt)^{\fl{\ap}{2}}$ has larger oscillations as $\ap$ increases.

Besides Lemma \ref{lemma1}, another important building block of our method is the properties of the Laplace operators as described below.
\begin{lemma}[{\bf Properties of the Laplace operators}]\label{lemma2}
For function $u$, assume the Laplacian function ${\mathcal U}(\bx) := (-\Dt)^{\fl{\ap}{2}}u(\bx)$ exists for $\bx \in {\mathbb R}^d$. 
Then it satisfies the following properties \cite{Samko}:
\bea\label{transit}
\displaystyle(-\Dt)^\fl{\ap}{2} \big[u(\bx - \bx_0)\big] = {\mathcal U}(\bx - \bx_0), \qquad \mbox{for}\ \ \ap \ge 0, 
\eea
for any point $\bx_0 \in {\mathbb R}^d$, and 
\bea\label{scaling}
\displaystyle (-\Dt)^\fl{\ap}{2}\big[u(\kappa \bx)\big] = |\kappa|^{\ap}{\mathcal U}(\kappa\bx), \qquad \mbox{for}\ \ \ap \ge 0, 
\eea
for  constant $\kappa \in {\mathbb R}$.
\end{lemma}
Lemma \ref{lemma2}  plays an important role in the design of our meshfree method, which allows us to find the analytical solution to the Laplacian of Gaussian RBFs with different shape parameters and center points. 
Combining (\ref{Eq-lemma1})--(\ref{scaling}), we immediately obtain that for any point $\bx_0 \in {\mathbb R}^d$ and shape parameter $\veps \in {\mathbb R}$, there is
\bea\label{core}
 (-\Dt)^{\fl{\ap}{2}} e^{-\veps^2|\bx-\bx_0|^2} = c_{d, \ap} |\veps|^\ap\,_1F_1\Big(\fl{d+\ap}{2}; \fl{d}{2}; -\veps^2|\bx - \bx_0|^2\Big),\qquad \mbox{for} \  \ \ap \ge 0.
 \eea
Here and in the following, we denote the constant $c_{d, \ap} = 2^\ap\Gamma((d+\ap)/2)/\Gamma(d/2)$, to be distinguished from $C_{d, \ap}$ in (\ref{integralFL}). 
It shows that over the entire space ${\mathbb R}^d$, the classical and fractional Laplacian of   Gaussian RBFs can be given in a unified form with parameter $\ap$.
As we will see in Section \ref{section2-2}, this uniform structure in (\ref{core}) provides the foundation to design the unified numerical methods for classical and fractional PDEs. 

\subsection{RBF discretization schemes}
\label{section2-2} 

In the following, we will present our meshfree method to approximate the Laplace operator $(-\Dt)^\fl{\ap}{2}$ and the related Poisson problems. 
Its generalization to time-dependent problems is straightforward (see Section \ref{section4-3}).
As mentioned previously, we will focus on the cases with $\ap \in (0, 2]$. 
Our method can be directly applied to discretize the operator $(-\Dt)^{2m}$ with $m \in {\mathbb N}$. 
Its generalization to the cases of $\ap > 2$ but $\ap \neq 2m$ requires the point-wise definition of  $(-\Dt)^\fl{\ap}{2}$, such as  (\ref{integralFL}) for $\ap \in (0, 2)$, which is beyond the scope of this work.

Let $\Og \subset {\mathbb R}^d$ be an open bounded domain. 
We consider the following Poisson problem  with Dirichlet boundary conditions: 
\bea
\label{BVP}
(-\Dt)^{\fl{\ap}{2}} u(\bx) = f(\bx),&&\mbox{for \ $\bx\in\Og$},\\
\label{BC}
 u(\bx) = g(\bx),&&\mbox{for \ $\bx\in \Upsilon$}. 
\eea
where we denote $\Upsilon = \p{\Og}$ for $\ap = 2$,\, or $\Upsilon = \Og^c = {\mathbb R}^d\backslash\Og$ if $\ap < 2$.  
If $\ap = 2$, the problem (\ref{BVP})--(\ref{BC})  becomes the classical Dirichlet Poisson problem.  
While $\ap \in (0, 2)$, it collapses to the fractional Poisson equation with extended Dirichlet boundary conditions on $\Og^c$. 
So far most studies on the fractional Poisson problem focus on the homogeneous Dirichlet boundary conditions (i.e., $g(\bx) \equiv 0$ in (\ref{BC})); see \cite{Duo2018, Duo-FDM2019, Acosta2017, Bonito2019} and references therein. 
Here, we consider more general boundary conditions $g(\bx)$. 
Usually,  if non-periodic boundary conditions are considered, the classical and fractional Poisson problems are discretized and solved separately. 
To the best of our knowledge, this is the first work to solve the classical and fractional Poisson problems in a single $\ap$-parametric scheme.  

Let $N$ and $\bar{N}$ be two positive integers, and $N < \bar{N}$. 
Denote $\bx_i$ (for $1 \le i \le \bar{N}$) as pre-defined collocation points on $\bar{\Og} = \Og \cup \p\Og$. 
For simplicity, we introduce
\beas
{\mathcal S}_\Og = \{\bx_i \in \Og \mid 1 \le i \le N\}, \qquad
{\mathcal S}_{\p\Og} = \{\bx_i \in \p{\Og} \mid N+1 \le i \le \bar{N}\}
\eeas
to represent the set of points in domain $\Og$ and on boundary $\p\Og$, respectively, and let  ${\mathcal S}_{\bar{\Og}} = {\mathcal S}_\Og\cup{\mathcal S}_{\p\Og}$. 
Assume that the function $u$ can be approximated by
\bea\label{Sol1D}
u(\bx) \approx \widehat{u}(\bx) :=\sum_{i=1}^{\bar{N}} \lambda_i\,\varphi^\veps(|\bx - \bx_i|), \qquad \mbox{for} \ \ \bx \in \bar{\Og}, 
\eea
where $\varphi^\veps(|\bx-\bx_i|)$ represents the Gaussian RBF  with shape parameter $\veps$ and center point $\bx_i$.  
For point $\bx \notin \bar{\Og}$,   we assume that $u$ satisfies (\ref{BC}). 
Then the coefficients $\lambda_i$ can be found by applying (\ref{Sol1D}) at a set of test points $\bx_k \in {\mathcal S}_{\bar{\Og}}$ that may or may not coincide with the center points. 
In the following, we will include a superscript to distinguish the center-point sets (i.e. ${\mathcal S}_\Og^c$ and ${\mathcal S}_{\p\Og}^c$) and test-point sets (i.e. ${\mathcal S}_\Og^t$ and ${\mathcal S}_{\p\Og}^t$).  
Note that  the number of test points should be the same as that of center points, i.e., $|\mathcal{S}_{\bar{\Og}}^t| = |{\mathcal S}_{\bar{\Og}}^c| = \bar{N}$, to ensure a square linear system for $\lambda_i$.

First, we will derive the approximation of the Dirichlet Laplacian, i.e., a finite dimensional representation of the operator $(-\Dt)^\fl{\ap}{2}$ with Dirichlet boundary conditions in (\ref{BC}).
To facilitate our explanation, we will start with separate discussion for $\ap = 2$ and $\ap < 2$,  demonstrating the difference between the classical and fractional Laplacians. 
Later, we will combine our results of $\ap = 2$ and $\ap < 2$ into a single $\ap$-parametric scheme. 
The situation of the classical Laplacian  is relatively simple, and combining (\ref{core}) with (\ref{Sol1D}) immediately leads to  the approximation: 
\bea\label{Eq-classical}
-\Dt_h  u(\bx) = -\Dt\widehat{u}(\bx) = c_{d, 2}|\varepsilon|^2\sum_{i = 1}^{\bar{N}}\lambda_i\, _1F_1\Big(\fl{d}{2}+1; \fl{d}{2}; -\varepsilon^2|\bx - \bx_i|^2\Big), \qquad \mbox{for} \ \ \bx \in\Og, 
\eea
where $-\Dt_h$ denotes the numerical approximation of the operator $-\Dt$. 

In contrast, the approximation to the fractional Laplacian (i.e. $\ap < 2$) is more complicated owing to its nonlocality over the entire space ${\mathbb R}^d$. 
For $\bx \in \Og$, we take the pointwise definition of the fractional Laplacian in (\ref{integralFL}) and reformulate it as 
\bea\label{Dlaplace}
(-\Dt)^{\fl{\ap}{2}}u(\bx) = C_{d, \ap} \bigg({\rm P.V.}\int_\Og\fl{{u}(\bx)- {u}(\by)}{|\bx - \by|^{d+\ap}} d\by + \int_{\Og^c}\fl{u(\bx) - u(\by)}{|\bx - \by|^{d+\ap}} d\by \bigg). \quad 
\eea
Substituting (\ref{Sol1D}) for $\bx \in \Og$ into (\ref{Dlaplace}) and taking the boundary conditions (\ref{BC}) into account, we obtain the approximation to the  fractional Laplacian as: 
\bea\label{Dlaplace1}
&&(-\Dt)^{\fl{\ap}{2}}_hu(\bx) = C_{d, \ap} \bigg({\rm P.V.}\int_\Og\fl{\widehat{u}(\bx)- \widehat{u}(\by)}{|\bx - \by|^{d+\ap}} d\by + \int_{\Og^c}\fl{\widehat{u}(\bx) - g(\by)}{|\bx - \by|^{d+\ap}} d\by \bigg) \nn\\
&&\hspace{2.05cm}= C_{d, \ap} \bigg({\rm P.V.}\int_{{\mathbb R}^d}\fl{\widehat{u}(\bx)- \widehat{u}(\by)}{|\bx - \by|^{d+\ap}} d\by + \int_{\Og^c}\fl{\widehat{u}(\by) - g(\by)}{|\bx - \by|^{d+\ap}} d\by \bigg) \qquad \qquad \nn\\
&&\hspace{2.05cm}= (-\Dt)^\fl{\ap}{2} \widehat{u}(\bx) +  C_{d, \ap}  \int_{\Og^c}\fl{\widehat{u}(\by) - g(\by)}{|\bx - \by|^{d+\ap}} d\by,
\eea
where the definition (\ref{integralFL}) is used again in the last line.  
We remark that substituting (\ref{Sol1D}) into (\ref{Dlaplace}) assumes as a default that $u(\bx) = \widehat{u}(\bx)$  for $\bx \in \Og^c$, however, the exact boundary condition is given by $u(\bx) = g(\bx)$ in (\ref{BC}).
Hence, the second term at the right side of (\ref{Dlaplace1}) can be viewed to match the difference of the Dirichlet boundary conditions on $\Og^c$, while the term $(-\Dt)^\fl{\ap}{2} \widehat{u}(\bx)$ in (\ref{Dlaplace1}) can be easily obtained by combining (\ref{core}) with (\ref{Sol1D}).

Combining (\ref{Eq-classical}) and (\ref{Dlaplace1}) yields a unified approximation to the Dirichlet Laplace operator $(-\Dt)^\fl{\ap}{2}$ for $\ap \in (0, 2]$, i.e. for $\bx \in \Og$, 
\bea\label{Eq-uniform}
(-\Dt)_h^{\fl{\ap}{2}}u(\bx) =  c_{d, \ap}|\varepsilon|^\ap\sum_{i = 1}^{\bar{N}}\lambda_i\, _1F_1\Big(\fl{d+\ap}{2}; \fl{d}{2}; -\varepsilon^2|\bx - \bx_i|^2\Big) +\zeta_\ap C_{d, \ap} \int_{\Og^c}\fl{\widehat{u}(\by) - g(\by)}{|\bx - \by|^{d+\ap}} d\by,\eea
where $\zeta_\ap = 1 - \lfloor\ap/2\rfloor$ with $\lfloor \cdot \rfloor$ being the floor function. 
Note that constants $c_{d, \ap}$ and $C_{d, \ap}$ are different and defined in (\ref{core}) and (\ref{integralFL}), respectively. 
The formulation in (\ref{Eq-uniform}) provides a uniform approximation to the classical and fractional Laplacian with the Dirichlet boundary conditions (\ref{BC}), and their difference lies in the integral term over $\Og^c$. 
If $\ap = 2$, the integral term in (\ref{Eq-uniform}) vanishes as $\lfloor \ap/2\rfloor = 1$, and (\ref{Eq-uniform}) reduces to the approximation of the classical Laplacian in (\ref{Eq-classical}). 
If $\ap < 2$, the nonlocal boundary conditions are {\it exactly} accounted through the integral over  $\Og^c$.
The approximation in (\ref{Eq-uniform}) again reveals the nonlocal nature of the fractional Laplacian -- the value at one point depends on all the other points over $\by \in {\mathbb R}^d$. 

Next, we will move to approximate the solution of the Poisson problem (\ref{BVP})--(\ref{BC}). 
Choose a set of test points $\bx_k \in {\mathcal S}_{\bar{\Og}}^t$.  
Substituting  (\ref{Eq-uniform}) with $\bx = \bx_k \in {\mathcal S}_\Og^t$ into the Poisson equation  (\ref{BVP}), we obtain the fully discretized scheme: 
\bea\label{Eq-gov}
&& \sum_{i = 1}^{\bar{N}} \lambda_i \bigg[{c}_{d, \ap}|\varepsilon|^\ap\,_1F_1\Big(\fl{d+\ap}{2}; \fl{d}{2};  -\veps^2|\bx_k - \bx_i|^2\Big) + \zeta_\ap C_{d, \ap}\int_{\Og^c}\fl{\varphi^\veps\big(|\by - \bx_i|\big)}{|\bx_k - \by|^{d+\ap}} d\by\bigg] \qquad\qquad \nn\\
&& \hspace{3cm} = f(\bx_k) + \zeta_\ap C_{d,\ap}\int_{\Og^c}\fl{g(\by)}{|\bx_k-\by|^{d+\ap}}d\by, \qquad \mbox{for} \ \ 1 \le k \le N. 
\eea
For $\bx_k \in {\mathcal S}_{\p\Og}^t$,  the boundary condition (\ref{BC}) leads to
\bea\label{Eq-BC}
\sum_{i  = 1}^{\bar{N}} \lambda_i \varphi^\veps\big(|\bx_k - \bx_i|\big)= g(\bx_k), \qquad 
\mbox{for} \ \ N+1 \le k \le  \bar{N}.
\eea
The discrete system (\ref{Eq-gov})--(\ref{Eq-BC}) has $\bar{N}$ equations with the same number of unknowns $\lambda_i$ (for $1\le i\le\bar{N}$).  
After obtaining $\lambda_i$, the solution of the Poisson problem (\ref{BVP})--(\ref{BC}) can be approximated from (\ref{Sol1D}). 
Even though nonlocal boundary conditions are imposed on $\Og^c$, the number of equations for the fractional Poisson problem remains the same as in the classical cases. 
For $\ap < 2$, the integrals over $\Og^c$ can be easily approximated, as their integrands  decay quickly and are free of singularities. 

Note that the linear system of (\ref{Eq-gov})--(\ref{Eq-BC}) has a full  stiffness matrix  for both classical and fractional Poisson problems, as the globally-supported Gaussian RBFs are used. 
Due to the nonlocal nature of the fractional Laplacian, evaluating integrals over $\Og^c$ and solving a linear system with full dense matrix are also required in other numerical methods \cite{Duo2018, Duo-FDM2019, Acosta2017};  thus our method does not introduce extra computations.
However,  compared to other local methods, our method can achieve higher accuracy with fewer number of points, implying fewer number of unknowns and smaller computational costs. 
This suggests that  the global numerical methods might be more beneficial for nonlocal or fractional problems.
\begin{remark}[{\bf Exact boundary conditions}] \label{remark1}
In the fractional cases, our method exactly incorporates the  nonlocal boundary condition (\ref{BC}) into the numerical scheme, which is one major difference from the method in \cite{Rosenfeld2019}. 
In \cite{Rosenfeld2019}, they  consider the boundary conditions on a small region  $\og\subset\Og^c$ and  assume  the Wendland RBF approximation to solution $u$ for all $\bx \in \Og \cup \og$. 
Consequently, their method not only introduces extra errors from boundary truncation but significantly increases the computational cost, as the fractional Poisson problem is actually solved on region $\Og \cup \og$ (in contrast to  $\Og \cup \p\Og$ in our method). 
Moreover, their method directly discretizes the pseudo-differential form of the fractional Laplacian in (\ref{pseudo}), reducing its usability on irregular domains.
\end{remark}

\section{Estimation of the Laplace operator}
\setcounter{equation}{0}
\label{section3}

In this section, we will study the performance of our method in approximating the Laplace operator $(-\Dt)^{\fl{\ap}{2}}$ for $\ap \in (0, 2]$. 
Let $\Og$ denote the domain of interest. 
{ Unless otherwise stated, we will choose the test points from the same set of center points, i.e., ${\mathcal S}_{\bar{\Og}}^t =  {\mathcal S}_{\bar{\Og}}^c$. 
This choice is not required by our method, but it has the advantage of yielding a  symmetric linear system of $\lambda_i$. } 
First, using the RBF approximation (\ref{Sol1D}) at all { test}  points $\bx_k  \in {\mathcal S}_{\bar{\Og}}^t$ gives a linear system of unknowns $\lambda_i$, where the coefficient matrix is positive definite with its entries given by the Gaussian RBFs $\varphi^\veps(|\bx_k - \bx_i|)$ for $\bx_i, \bx_k  \in {\mathcal S}_{\bar{\Og}}^c$. 
Solving and substituting  $\lambda_i$ into (\ref{Eq-uniform}), we then obtain the approximation of $(-\Dt)^\fl{\ap}{2}u$ for any  $\bx \in \Og$. 
For $\ap < 2$, the integral term in  (\ref{Eq-uniform}) can be computed by setting  $g(\bx) = u(\bx)$ for $\bx \in \Og^c$. 

In the following, we will study and compare numerical errors under different conditions of $u$, where numerical errors are computed as the root mean square (RMS) error, i.e., 
\beas\label{rms-operator}
\|e_\Dt\|_{\rm rms} &=& \bigg(\fl{1}{M}\sum_{j =1}^M \Big[(-\Dt)^{\fl{\ap}{2}}u(\bx_j) - (-\Dt)_h^{\fl{\ap}{2}}u(\bx_j)\Big]^2\bigg)^{1/2}
\eeas
with $\bx_j$ (for $1\le j \le M$)  denoting the interpolation points on domain $\Og$. 
To better estimate numerical errors,  the number of interpolation points $M$ is chosen to be much larger than that of center points, i.e., $M \gg \bar{N}$. 
We usually take a large enough $M$ such that the error $\|e_\Dt\|_{\rm rms}$ is insensitive to the number of interpolation points.

\subsection{Globally smooth functions}
\label{section3-1}

Consider a globally smooth  function $u(x) = 1/(1+x^2)$ for $x \in {\mathbb R}$. 
In this case, the function $(-\Dt)^{\fl{\ap}{2}}u$ can be exactly given by: 
\bea\label{fun3-1}
(-\Dt)^\fl{\ap}{2}u(x)= \Gamma\big(1+\ap\big)\,_2F_1\Big(\fl{1+\ap}{2}, \, \fl{2+\ap}{2};\, \fl{1}{2}; \, -x^2\Big), \qquad \mbox{for}\ \ \ap \ge 0, \ \ x \in{\mathbb R},
\eea
where $_2F_1$ represents the Gauss hypergeometric  function. 
We remark that  the exact solution (\ref{fun3-1}) holds for any $\ap \ge 0$. 
It is easy to verify that the result (\ref{fun3-1}) is consistent with the classical integer-order derivatives if $\ap = 2m$ for $m \in {\mathbb N}$.
\begin{figure}[htb!]
\centerline{\includegraphics[height = 5.76cm, width = 7.86cm]{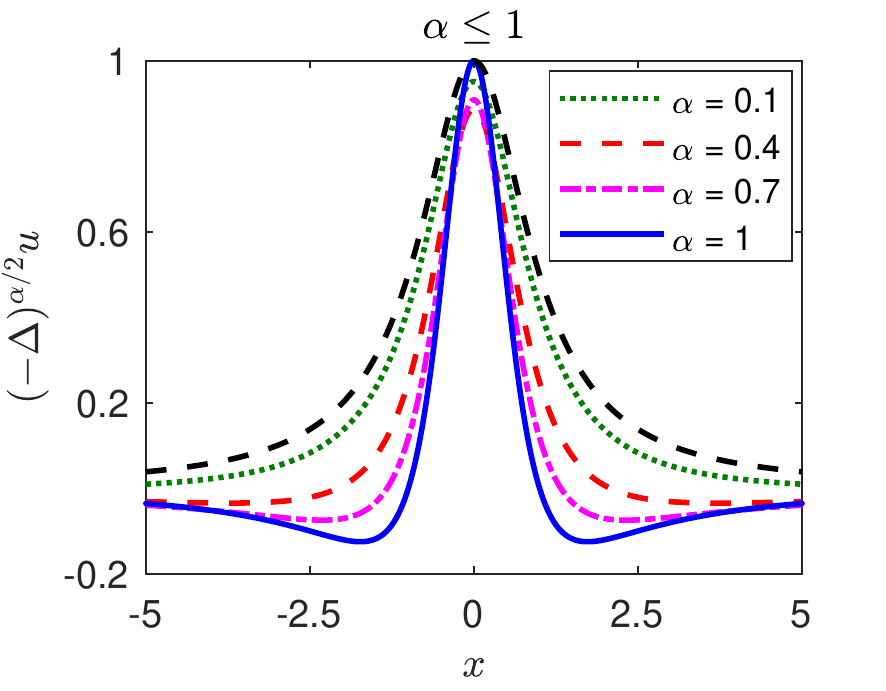}\hspace{-3mm}
\includegraphics[height = 5.76cm, width = 7.86cm]{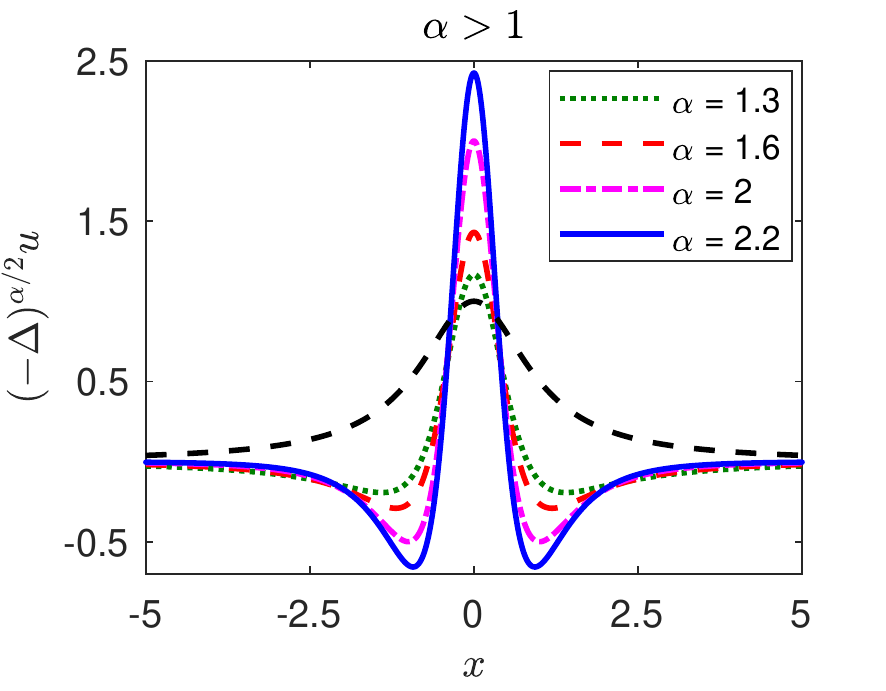}}
\caption{Illustration of  function $(-\Dt)^\fl{\ap}{2}(1+x^2)^{-1}$ for various $\ap$. 
For easy comparison, the function $u(x) = (1+x^2)^{-1}$ is also presented by  black dashed lines.} \label{Fig3-1-0}
\end{figure}
Fig. \ref{Fig3-1-0} illustrates  function $(-\Dt)^\fl{\ap}{2}u$  for various $\ap$. 
It shows that  function $(-\Dt)^\fl{\ap}{2}u$ approaches to $u$ as $\ap \to 0$. 
The nonlocal effects become stronger for smaller $\ap$, and particularly we find that $(-\Dt)^\fl{\ap}{2}u(x) \le u(x)$ for any $x \in {\mathbb R}$ if $\ap \le 1$.
For $\ap > 2$ but $\fl{\ap}{2} \notin {\mathbb N}$, the operator $(-\Dt)^\fl{\ap}{2}$ remains nonlocal, which is beyond the scope of our study. 

Table \ref{Tab3-1} presents the RMS errors of our method in approximating $(-\Dt)^{\fl{\ap}{2}}u$ on domain $\Og = (-2, 2)$. 
Here,  we take the shape  parameter $\varepsilon = 2$.  
The center and test points are chosen to be uniformly distributed on $[-2, 2]$. 
{ We remark that our method is flexible in choosing center and test points, but the ``best" shape parameter might change according to this choice. 
For example, our studies show that using the Chebyshev points as RBF center and test points can give the similar numerical accuracy as in Table \ref{Tab3-1},  if a larger shape parameter is used. }
Note that  the optimal choice of shape parameter and center and test points of RBF-based methods is still an open research topic \cite{Majdisova2017, Fasshauer2007}, and we will leave it for our future study \cite{Wu0020}.
\begin{table}[htb!]
\begin{center} 
\begin{tabular}{|c||c|c|c|c|c|}
\hline
$\bar{N}$ & $\ap  = 0.4 $ & $\ap = 1$ & $\ap = 1.6$ & $ \ap  = 2$ & $\mathcal{K}$ \\
\hline
9 & 1.957E-3 & 2.177E-2& 8.091E-2& 1.941E-1 & 5.431\\ 
\hline
17 & 8.442E-4 & 4.009E-3& 2.230E-2& 8.116E-2 & 5.079E3\\  
\hline
33 & 1.010E-6& 7.856E-6& 7.732E-5& 4.949E-4 & 2.369E14\\ 
\hline
65 & 2.220E-9& 1.486E-8& 1.832E-7& 1.514E-6 & 2.350E17 \\
\hline
\end{tabular}
\caption{Numerical errors $\|e_\Dt\|_{\rm rms}$ and condition number ${\mathcal K}$ in approximating $(-\Dt)^\fl{\ap}{2}u$ for $x \in (-2, 2)$, where the shape parameter $\varepsilon = 2$ and the exact solution is given in (\ref{fun3-1}).} \label{Tab3-1}
\end{center}
\end{table}

Table \ref{Tab3-1}  shows that our method yields a good approximation to  $(-\Dt)^{\fl{\ap}{2}}u$ even with a small number of points $\bar{N}$.
Comparing the errors of different $\ap$, we find that the larger the exponent $\ap$, the bigger the numerical errors, { but a spectral accuracy is achieved for any $\ap \in (0, 2]$. 
The condition number ${\mathcal K}$ of the linear system increases with the number of points $\bar{N}$, which may lead to an ill-conditioned system if $\bar{N}$ is too big. 
The ill-conditioning is one issue of  methods with infinitely differentiable RBFs, and so far different strategies have been developed to improve or control it (see e.g. \cite{Kansa2017} and references therein). 
Recently,  new algorithms based on expansion of RBFs have been also proposed in \cite{Fasshauer2012, Fornberg2011, Kormann2019} to tackle the ill-conditioning issues. 
However, the expansion of Gaussian RBFs destroys the  property in Lemma \ref{lemma1} and thus fails to work for our method.
In practical simulations, we can control the condition number  to  ${\mathcal O}(10^{13}) \sim {\mathcal O}(10^{17})$ via adjusting the shape parameter $\veps$ so as to obtain the best accuracy. 
Furthermore, multi-precision toolboxes and domain decompositions are also recommended in the literature \cite{Sarra2011, Kansa2017, Sarra2017}, if higher accuracy is demanded.}

\begin{figure}[htb!]
\centerline{(a)\includegraphics[height = 5.76cm, width = 7.86cm]{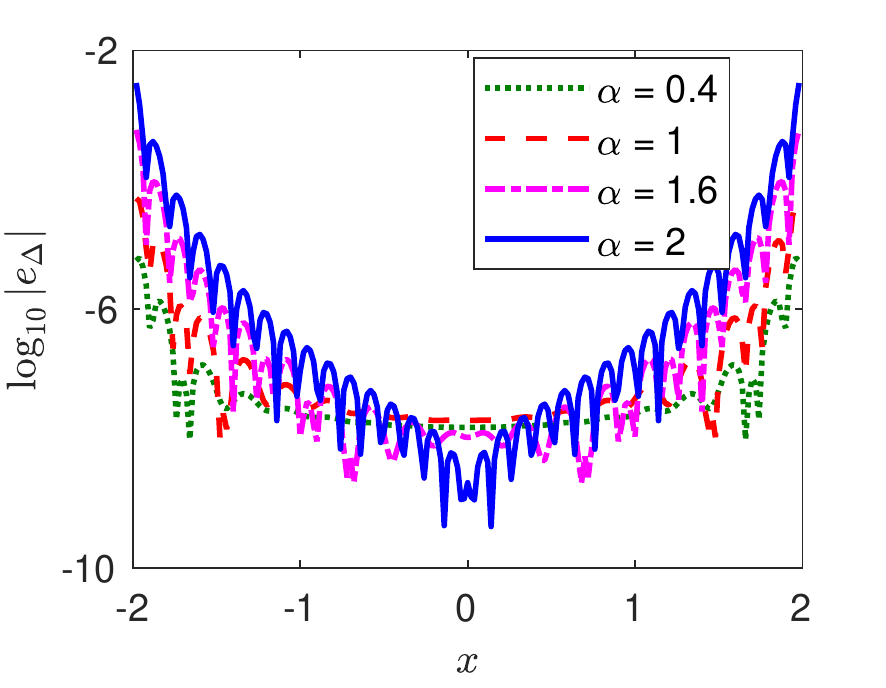}\hspace{-5mm}
(b)\includegraphics[height = 5.76cm, width = 7.86cm]{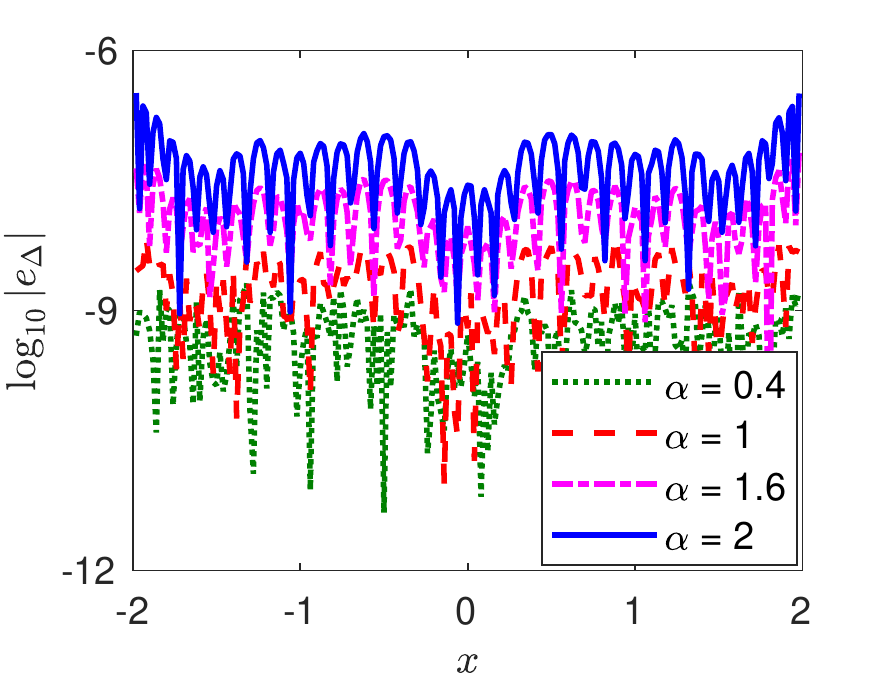}}
\caption{Error distribution in approximating (\ref{fun3-1}), where $|e_\Dt| = \big|(-\Dt)^\fl{\ap}{2}u - (-\Dt)_h^\fl{\ap}{2}u\big|$. (a) $\bar{N} = 33$ and (b) $\bar{N} = 129$. } \label{Fig3-1-1}
\end{figure}

As discussed previously, our method of approximating the classical and fractional Laplacians are the same, where extra efforts  are required in the fractional cases to evaluate the integrals over the domain ${\mathbb R}\backslash(-2, 2)$. 
In Fig. \ref{Fig3-1-1}, we further demonstrate the pointwise error  for $\bar{N}  = 33$ and $129$. 
It shows that  { if $\bar{N}$ is small} the maximum error occurs symmetrically around the domain boundary (see Fig. \ref{Fig3-1-1} (a)). 
{ This is simply because of the lack of points around the domain boundary.} 
Our extensive studies show that including more RBF points around or outside of the boundary could improve the accuracy of approximation, { consistent with the observations in \cite{Fedoseyev2002}. }

\begin{figure}[htb!]
\centerline{
(a)\hspace{-0.5mm}\includegraphics[height = 3.86cm, width = 4.46cm]{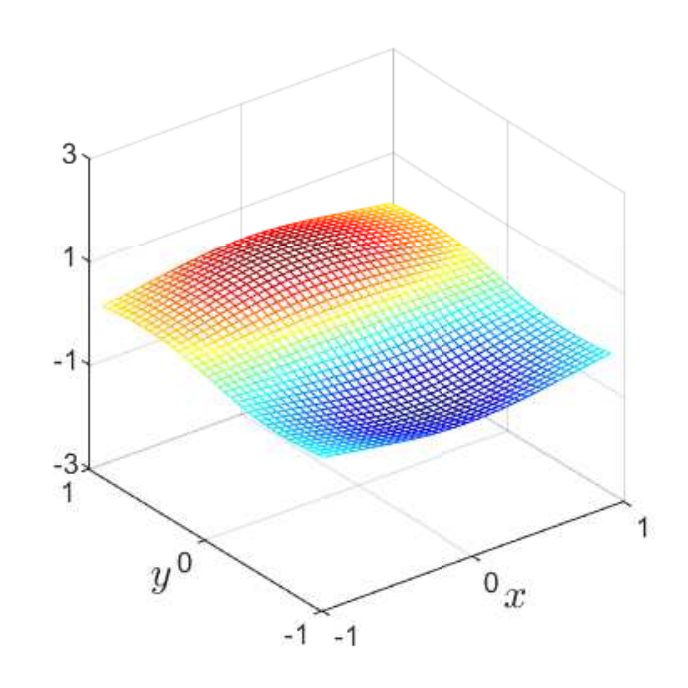}
\hspace{-2mm}
(b)\hspace{-0.5mm}\includegraphics[height = 3.86cm, width = 4.46cm]{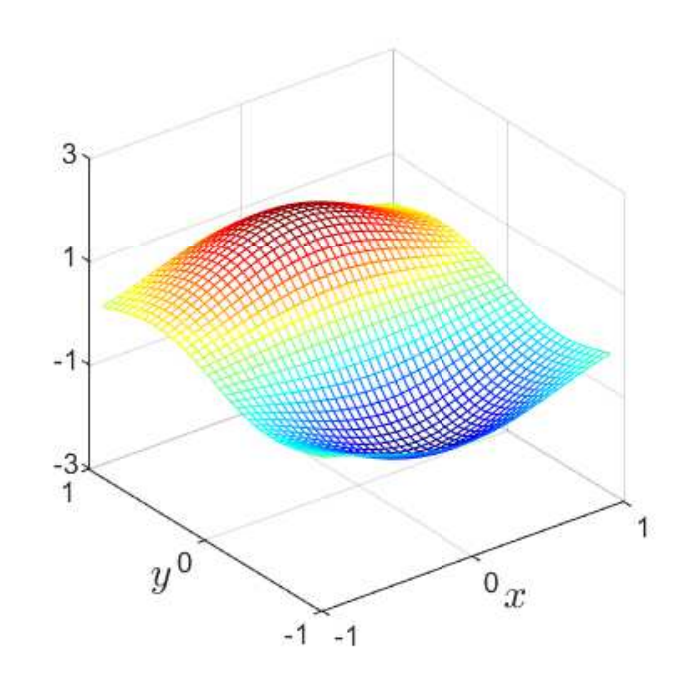}
\hspace{-2mm}
(c)\hspace{-0.5mm}\includegraphics[height = 3.86cm, width = 4.46cm]{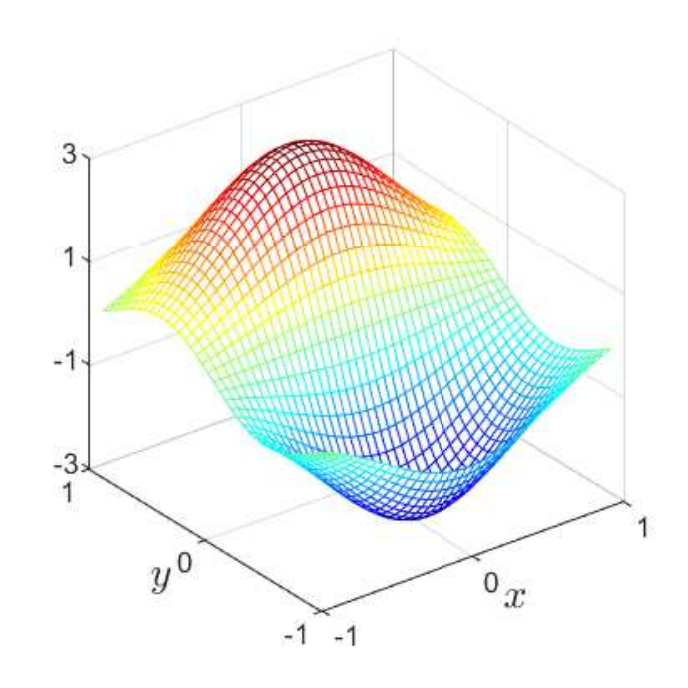}}
\caption{Numerical approximation of $(-\Dt)^{\fl{\ap}{2}}u$ for two-dimensional function $u(\bx) = e^{-(x^2+y^2)}\sin(y)$.  From (a) to (c):  $\ap = 0.6, 1.4$, and $2$. }\label{Fig3-1-2}
\end{figure}

In Fig. \ref{Fig3-1-2}, we present the numerical approximation of the two-dimensional function $(-\Dt)^{\fl{\ap}{2}}u$ with $u(\bx) = \exp[-(x^2+y^2)]\sin(y)$. 
It shows that the result is symmetric with respect to the $x$-axis, as the function $u$.  
In this case, the exact solution of $(-\Dt)^{\fl{\ap}{2}}u$ is unknown for $\ap < 2$, but if  $\ap = 2$ we have the analytical result: 
\beas
-\Dt u(\bx) = e^{-(x^2+y^2)} \Big[\sin(y)\big(4x^2+4y^2-5\big) -4y\cos(y)\Big], \qquad\mbox{for} \ \bx \in {\mathbb R}^2. 
\eeas
For $\ap = 2$, our approximate results agree well with their exact solutions (see Fig. \ref{Fig3-1-2} (c)). { Furthermore, we compare our method with the finite difference method in Table \ref{table-extra}, where the RBF center points are taken to be the same as the finite difference grid points.
It shows that  our method provides a more accurate approximation with the same number of points $\bar{N}$. 
Moreover,  the geometric flexibility and easy implementation make our method more advantageous in high dimensions. }
\begin{table}[h!]
\begin{center}
\begin{tabular}{|ccccc|}
\hline
$\bar{N}$  & $4^2$ & $5^2$ & $6^2$ & $7^2$\\
\hline
FDM       & 1.046E-2  & 5.493E-3 & 3.388E-3 &2.299E-3\\
\hline
RBF        & 3.094E-3 & 5.479E-5 & 5.818E-7 &4.010E-9 \\
\hline
\end{tabular}
\caption{Numerical errors $\|e_\Dt\|_{\rm rms}$ of the finite difference method (FDM) and our method (RBF) with  $\veps = 1$  in approximating function $-\Dt u$ on $(-1, 1)^2$, where $u(\bx) = e^{-(x^2+y^2)}\sin(y)$. }
\label{table-extra}
\end{center}
\end{table}

\subsection{Compactly supported functions}
\label{section3-2}

In the following, we approximate the Laplacian of compactly supported functions which are often studied in the field of fractional calculus. 
Consider  function $u(x) = x(1-x^2)^p_+$ for $x\in{\mathbb R}$, which has compact support on $(-1, 1)$ and is always zero for $x \notin (-1, 1)$. 
For $p > -1$, there is the exact solution \cite{Dyda2012}:  
\bea\label{exact-compact}
(-\Dt)^{\fl{\ap}{2}} u(x) = \fl{2^\ap (\ap+1)\Gamma((1+\ap)/2)\Gamma(p+1)}{\sqrt{\pi}\,\Gamma(p+1-\ap/2)}\,_2F_1\Big(\fl{\ap+3}{2},\, -p+\fl{\ap}{2}; \, \fl{3}{2}; \, x^2\Big)x, \quad\mbox{for \ $|x| <1$}.
\eea  
It shows in  \cite{Dyda2012} that the exact solution (\ref{exact-compact}) holds for $\ap \in (0, 2)$, but we find that  it is also valid for $\ap = 2m$ with $m \in {\mathbb N}$. 
\begin{figure}[htb!]
\centerline{\includegraphics[height = 5.76cm, width = 7.86cm]{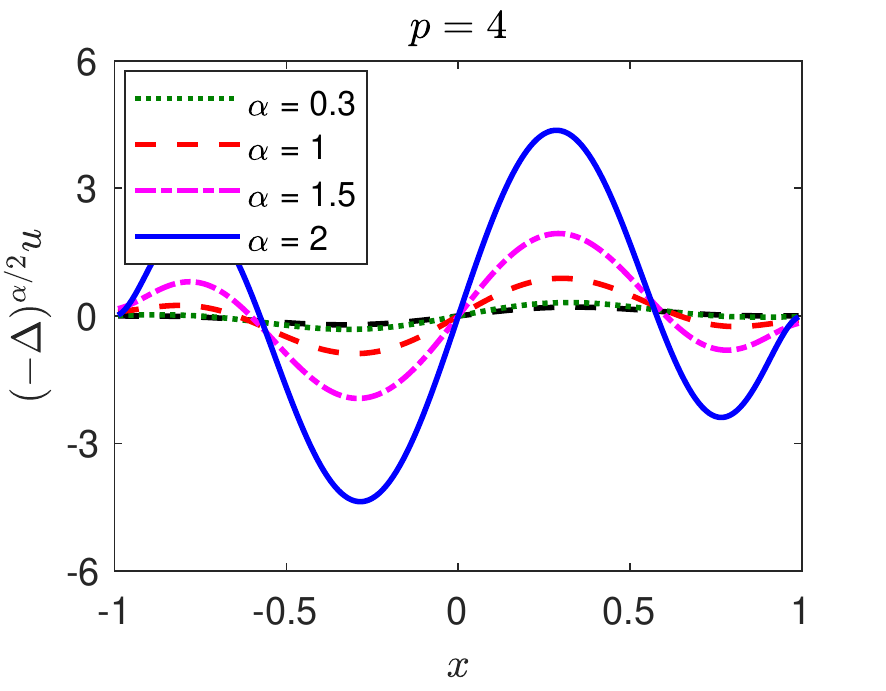}
\hspace{-3mm}
\includegraphics[height = 5.76cm, width = 7.86cm]{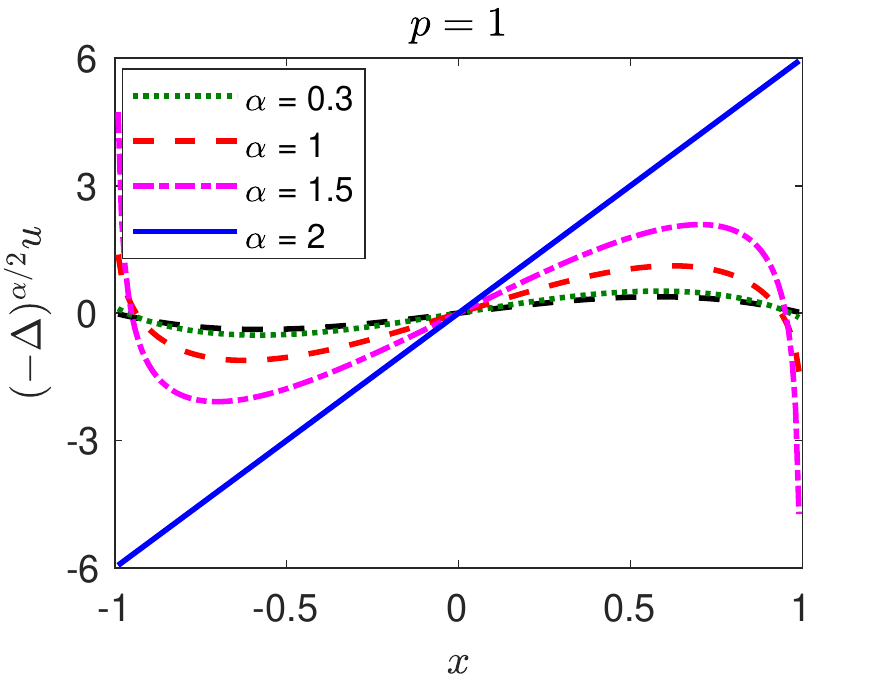}}
\caption{Illustration of function $(-\Dt)^\fl{\ap}{2}\big[x(1-x^2)^{p}_+\big]$ for different $p$ and $\ap$. 
For easy comparison, the function $u(x) = x(1-x^2)^{p}_+$ is also presented by black dashed lines.}\label{Fig3-2-0}
\end{figure}
It is easy to see that  the { differentiability of function $u$ at $x = \pm 1$ increases as $p > 0$ increases.}
Fig. \ref{Fig3-2-0} compares the solution behavior of $(-\Dt)^\fl{\ap}{2}\big[x(1-x^2)^{p}_+\big]$ for different $p$ and $\ap$. 
It shows that for $p = 1$,  the results become much sharper around the boundary as $\ap < 2$ increases, which makes the accurate approximation more challenging { for low-order methods}.

\begin{figure}[htb!]
\centerline{
(a)\includegraphics[height = 5.76cm, width = 7.86cm]{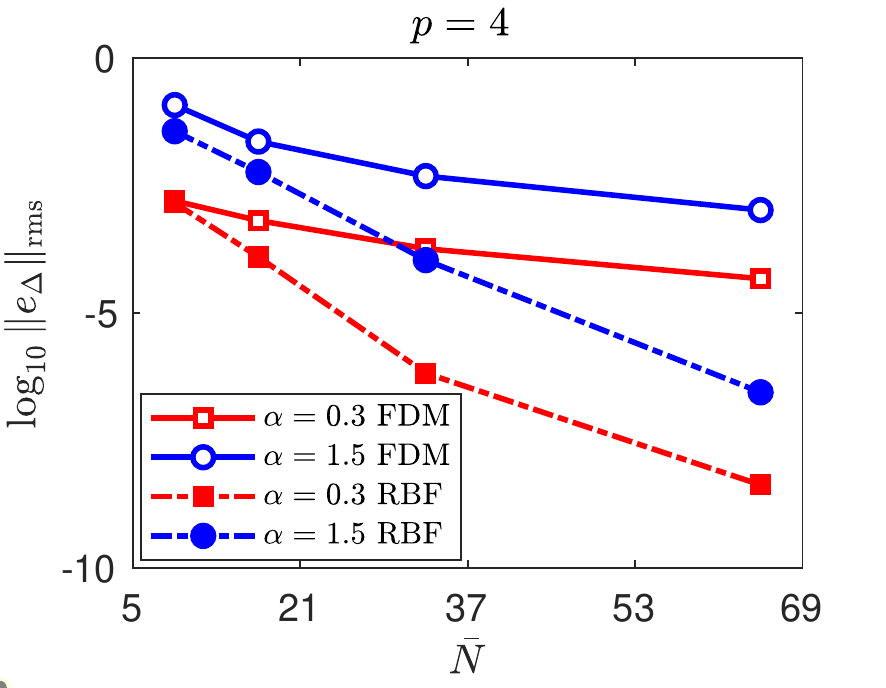}
\hspace{-5mm}
(b)\includegraphics[height = 5.76cm, width = 7.86cm]{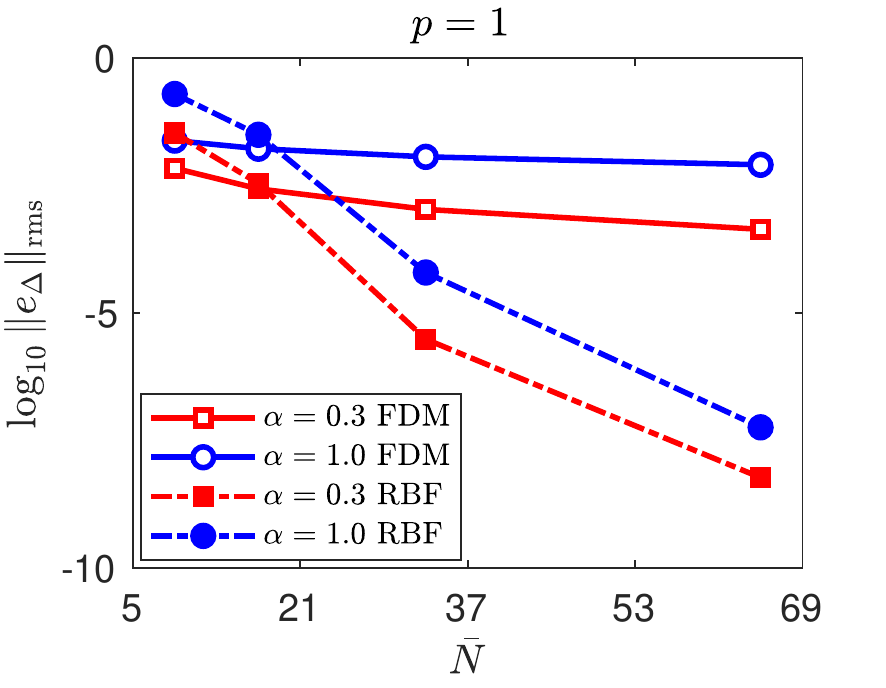}}
\caption{Comparison of numerical errors from our method with $\veps = 4$ and the finite difference method in approximating $(-\Dt)^\fl{\ap}{2}\big[x(1-x^2)^{p}_+\big]$ for $x \in (-1, 1)$. }\label{Fig3-2-2-extra}
\end{figure}

Table \ref{Tab3-2-1} presents the RMS errors for cases of $p = 4$ and $p = 1$, where the RBF center and test points are chosen uniformly on  $[-1, 1]$ and the shape parameter $\varepsilon = 4$.
The function $(-\Dt)^\fl{\ap}{2}u$ is approximated following the same process as in Section \ref{section3-1}, except the boundary conditions considered here { are} $g(x) = 0$ for $x\in {\mathbb R}\backslash(-1,1)$.
It shows that our method provides a good approximation to $(-\Dt)^\fl{\ap}{2}u$ for both $\ap  = 2$ and $\ap < 2$. 
\begin{table}[htb!]
\begin{center} 
\begin{tabular}{|c|c|c|c||c||c|c|c|c|}
\hline
\multicolumn{4}{|c||}{$p  = 4$} & &\multicolumn{4}{c|}{$p  = 1$}\\
\cline{1-4}
\cline{6-9}
$\ap = 0.3$ & $\ap = 1$ & $\ap = 1.5$ & $\ap  = 2$ & $\bar{N}$ & $\ap = 0.3$ & $\ap = 1$ & $\ap = 1.5$ & $\ap  = 2$\\
\hline
1.208E-1 & 4.521E-1 & 1.2041     & 3.2937     & 5   & 1.725E-1 & 6.789E-1 & 2.0617    & 6.5544\\
\hline
1.456E-3 & 9.267E-3 & 3.709E-2 & 1.539E-1 & 9   & 3.423E-2 & 1.993E-1 & 7.551E-1 & 2.9834\\
\hline
1.266E-4 & 1.169E-3 & 5.867E-3 & 2.963E-2 & 17 & 3.519E-3 & 3.159E-2 & 1.551E-1 & 7.596E-1\\
\hline
6.636E-7 & 1.389E-5 & 1.096E-4 & 8.462E-4 & 33 & 3.074E-6 & 6.261E-5 & 4.905E-4 & 3.777E-3\\
\hline
4.251E-9 & 4.174E-8 & 2.793E-7 & 2.029E-6 & 65 & 5.917E-9 & 5.736E-8 & 3.803E-7 & 2.720E-6\\
\hline
\end{tabular}
\caption{Numerical errors $\|e_\Dt\|_{\rm rms}$ in approximating $(-\Dt)^{\fl{\ap}{2}}u$ for $x \in (-1, 1)$, where  the shape parameter $\varepsilon = 4$ and the exact solution is given in 
(\ref{exact-compact}).} \label{Tab3-2-1}
\end{center}
\end{table}
The larger the exponent $\ap$, the bigger the numerical errors, similar observation to that from Table \ref{Tab3-1}. 
Our studies show that { the} numerical errors are symmetric with respect to $x = 0$. Compared to the finite difference methods in \cite{Duo2018, Duo-FDM2019}, the proposed method has significantly smaller errors with the same { number of points (see Fig. \ref{Fig3-2-2-extra}}). 
{ Moreover, the finite difference method fails to converge when $p = 1$ and $1 < \ap < 2$, since the function $u$ in this case does not satisfy the consistency condition as discussed in \cite{Duo-FDM2019, Duo-TFL2019}. 
In contrast to it, our method achieves a spectral accuracy for both $p = 1$ and $4$.}

\begin{figure}[htb!]
\centerline{
(a)\includegraphics[height = 5.76cm, width = 7.86cm]{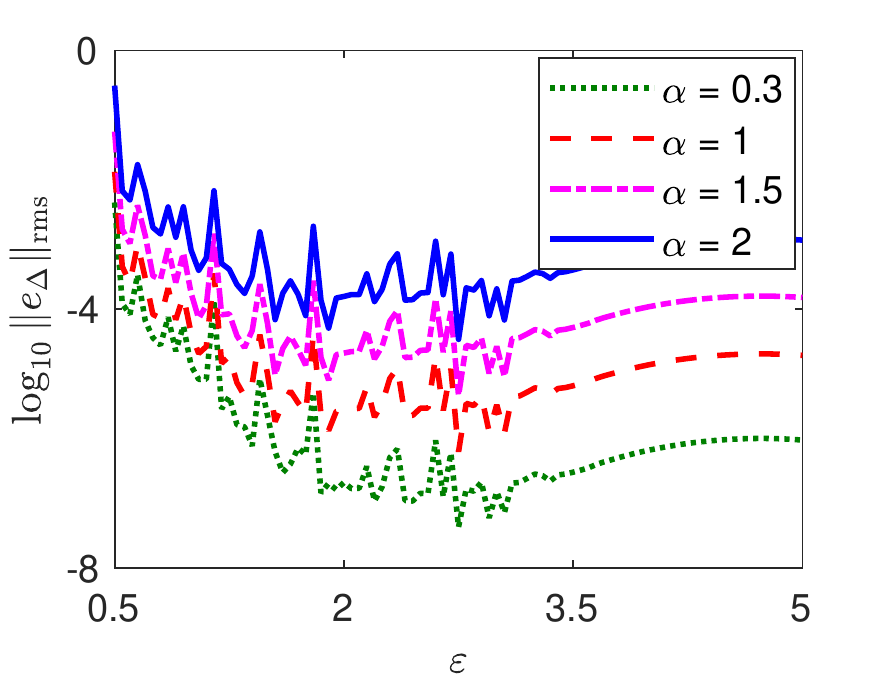}
\hspace{-5mm}
(b)\includegraphics[height = 5.76cm, width = 7.86cm]{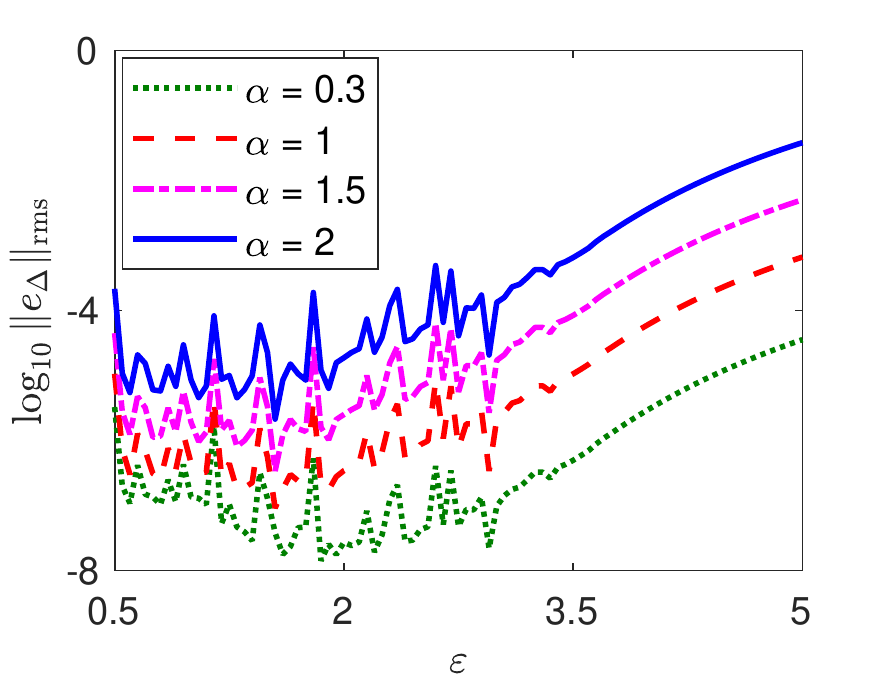}}
\caption{Numerical errors versus shape parameters in approximating the function in (\ref{exact-compact}), where the number of points $\bar{N} = 33$. (a) $p = 4$; (b) $p = 1$.}\label{Fig3-2-2}
\end{figure}

In Fig. \ref{Fig3-2-2}, we further study the numerical errors for different shape parameters, where the number of points $\bar{N} = 33$ is fixed. 
We find that the dependence of numerical errors on the shape parameter are qualitatively the same for different $\ap$ and $p$.  
To obtain good accuracy, one should choose the shape parameter neither too small nor too large.  
The optimal shape parameter might exist and depend not only on the choice RBF points but also on  function $u$ and exponent $\ap$. 
In Fig.  \ref{Fig3-2-2}, for example, the optimal shape parameter occurs around $2.8$ for $p  = 4$, while $1.8$ for $p = 1$. 
How to find the optimal shape parameter remains an active research topic in the field of RBF-related methods, and we will leave it for our future research. 
\section{Solutions of classical and fractional PDEs}
\label{section4}
\setcounter{equation}{0}

In this section, we test the performance of our method in solving both classical and fractional PDEs. 
So far, most existing numerical methods for  fractional PDEs with the integral fractional Laplacian are incompatible with those for their classical counterparts, due to different formulations of the Laplace operators as well as the boundary conditions. 
For example,  the weak formulation of the classical and fractional Poisson equations are significantly different when using finite element methods. 
Consequently, numerical schemes and computer codes have to be developed separately for the classical and fractional problems.  
One  important merit of our method is to solve both the classical and fractional problems in a single scheme, owing to the unified Laplace formulation of the Gaussian RBFs in Lemmas \ref{lemma1} and \ref{lemma2}. 
We will further demonstrate this and other advantages of our method in this section. 
To quantify its performance, we define { the} RMS error in solution $u$ as: 
\beas
\|e_u\|_{\rm rms} = \bigg(\fl{1}{M}\sum_{j = 1}^M \big|u_j - u_j^h\big|^2\bigg)^{1/2}, 
\eeas 
where $M$ denotes the number of interpolation points on $\Og$,  and $u_j$ and $u_j^h$,  respectively, represent the exact and numerical solutions at point  $\bx_j \in {\Og}$.

\subsection{One-dimensional Poisson problems}
\label{section4-1}

In this case,  we take $d = 1$ in (\ref{BVP})--(\ref{BC}) and choose the domain $\Og = (-1, 1)$.   
We study a benchmark fractional Poisson problem in \cite{Duo2018, Rosenfeld2019, Acosta2017}, i.e., choosing 
\bea\label{fun4-1}
f(x) = \fl{2^\ap\Gamma(\fl{\ap+1}{2})\Gamma(s+1+\fl{\ap}{2})}{\sqrt{\pi}\Gamma(s+1)} \,_2F_1\Big(\fl{\ap+1}{2}, -s;\,  \fl{1}{2}; \, x^2\Big), &\quad& \mbox{for} \  x \in \Og, \qquad \\
\label{fun4-2}
g(x) \equiv 0,  &\quad& \mbox{for} \  x \in \Upsilon,
\eea
for $s \in {\mathbb N}^0$. 
Noticing the definition of $\Upsilon$, (\ref{fun4-2}) implies that the two-point zero boundary conditions at $x = \pm1$ are imposed for the classical problem with $\ap = 2$, while the extended homogeneous Dirichlet boundary conditions are considered for $\ap < 2$. 
The exact solution of the Poisson equation with (\ref{fun4-1})--(\ref{fun4-2})  is given by $u(x) = (1-x^2)_+^{s+\fl{\ap}{2}}$ for any $\ap \in (0, 2]$. { It is evident that} the larger the value of $s$, the better the regularity of solution $u$ up to the boundary.

\begin{figure}[htb!]
\centerline{(a)\includegraphics[height = 5.76cm, width = 7.86cm]{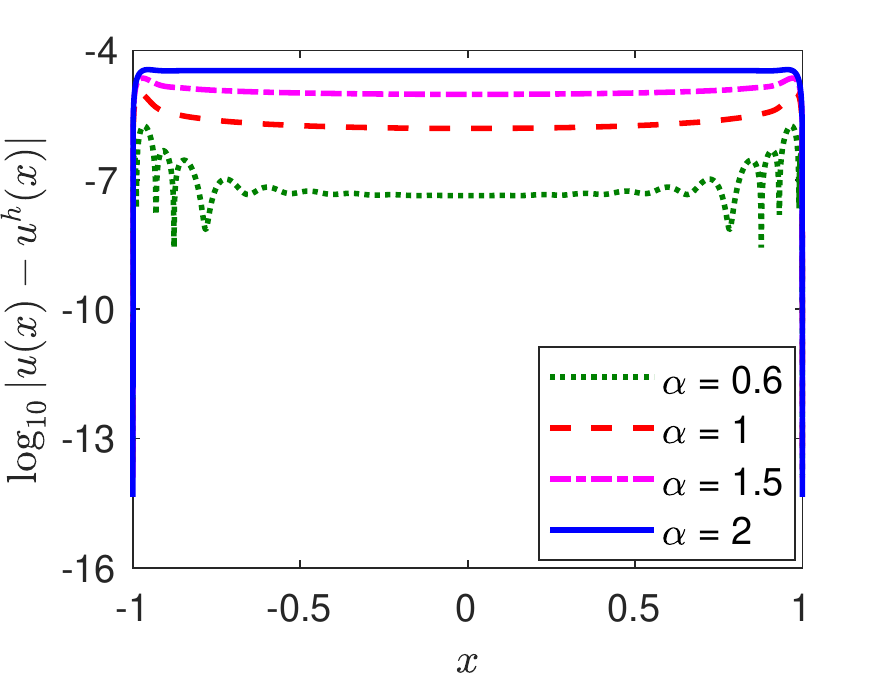}\hspace{-3mm}
(b)\includegraphics[height = 5.76cm, width = 7.86cm]{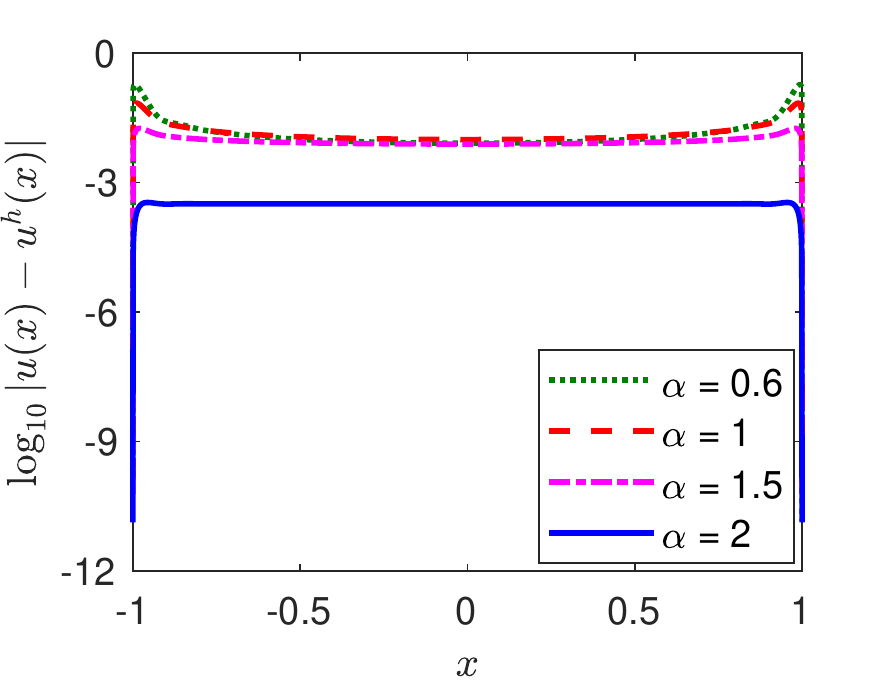}}
\caption{Numerical errors in solving the 1D Poisson problem with $f$ and $g$ defined in (\ref{fun4-1})--(\ref{fun4-2}), where we use $\varepsilon = 4.5$ and $\bar{N} = 33$ in our simulations. (a) $s  = 3$; (b) $s  = 0$.}\label{Fig4-1-1}
\end{figure}

Table \ref{Tab4-1-1} presents the RMS errors $\|e_u\|_{\rm rms}$ when different $s$ are chosen in (\ref{fun4-1}). 
In our simulations, the center and test points are taken to be  uniformly distributed on $[-1, 1]$, and the shape parameter  $\varepsilon = 4.5$ is used.
From Table \ref{Tab4-1-1}, we find that numerical errors depend on the solution regularity and exponent $\ap$. 
\begin{table}[htb!]
\begin{center} 
\begin{tabular}{|c|c|c|c||c||c|c|c|c|}
\hline
\multicolumn{4}{|c||}{$s  = 3$} & &\multicolumn{4}{c|}{$s  = 0$}\\
\cline{1-4}
\cline{6-9}
$\ap = 0.6$ & $\ap = 1$ & $\ap = 1.5$ & $\ap  = 2$ & $\bar{N}$ & $\ap = 0.6$ & $\ap = 1$ & $\ap = 1.5$ & $\ap  = 2$\\
\hline
3.262E-1& 4.134E-1 &4.791E-1 & 5.074E-1 & $5$ & 5.863E-1 & 6.539E-1 & 6.922E-1 & 6.952E-1 \\
\hline
5.052E-3& 1.308E-2 &4.263E-2 & 1.166E-1 &$9$  & 1.627E-1 & 1.357E-1 & 1.437E-1 & 2.098E-1 \\
\hline
1.147E-4& 1.648E-4 &2.553E-4 & 3.750E-4 &$17$  & 7.746E-2 & 5.225E-2 & 3.405E-2 & 1.607E-2 \\
\hline
3.120E-7& 2.909E-6 &1.174E-5 & 2.383E-5 &$33$ & 3.265E-2 & 1.867E-2 & 9.485E-3  & 3.193E-4 \\
\hline
8.147E-8& 1.719E-7 & 2.288E-7 & 1.462E-6 &$65$ & 1.631E-2 & 8.538E-3 & 3.966E-3 & 3.014E-6\\
\hline
\end{tabular}
\caption{Numerical errors $\|e_u\|_{\rm rms}$ in solving the 1D Poisson problem with $f$ and $g$ defined in (\ref{fun4-1})--(\ref{fun4-2}), where  the shape parameter $\varepsilon = 4.5$.}\label{Tab4-1-1}
\end{center}
\end{table}
If the solution $u$ is smooth enough up to the boundary (e.g., $s = 3$), the numerical errors decrease { with a spectral rate} as the number of RBF points increases. 
Moreover, our numerical errors are much smaller than those computed from finite difference method \cite[Tables 4--5]{Duo2018} with the same number of points. 
While $s = 0$, the solution of the fractional Poisson problem has low regularity around the boundary, which negatively affects the accuracy of our method. 
Consequently, the numerical errors of $s = 0$ are much larger  than those of $s = 3$ for any $\ap \in (0, 2)$. 
However, the effect of $s$ on the solution of the classical problem (i.e., $\ap = 2$) is less significant, as the classical Laplacian is a local operator.  
For both $s = 0$ and $3$, the numerical errors are generally larger around the two boundary points $x = \pm 1$ (see Fig.  \ref{Fig4-1-1}). 
It implies that including more points around the boundary might improve the accuracy of our method. 

\begin{figure}[htb!]
\centerline{(a)\includegraphics[height = 5.76cm, width = 7.86cm]{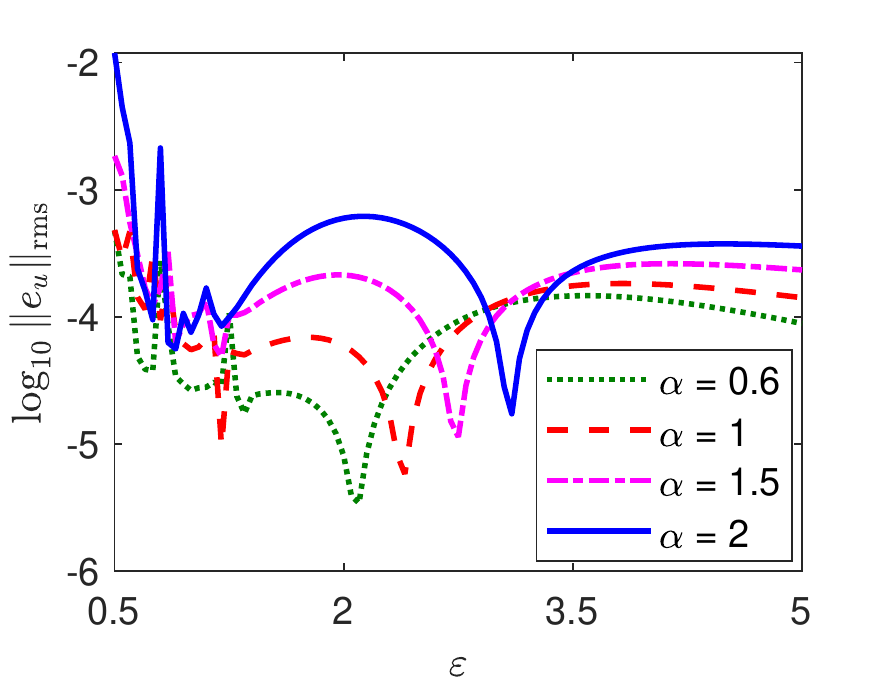}\hspace{-3mm}
(b)\includegraphics[height = 5.76cm, width = 7.86cm]{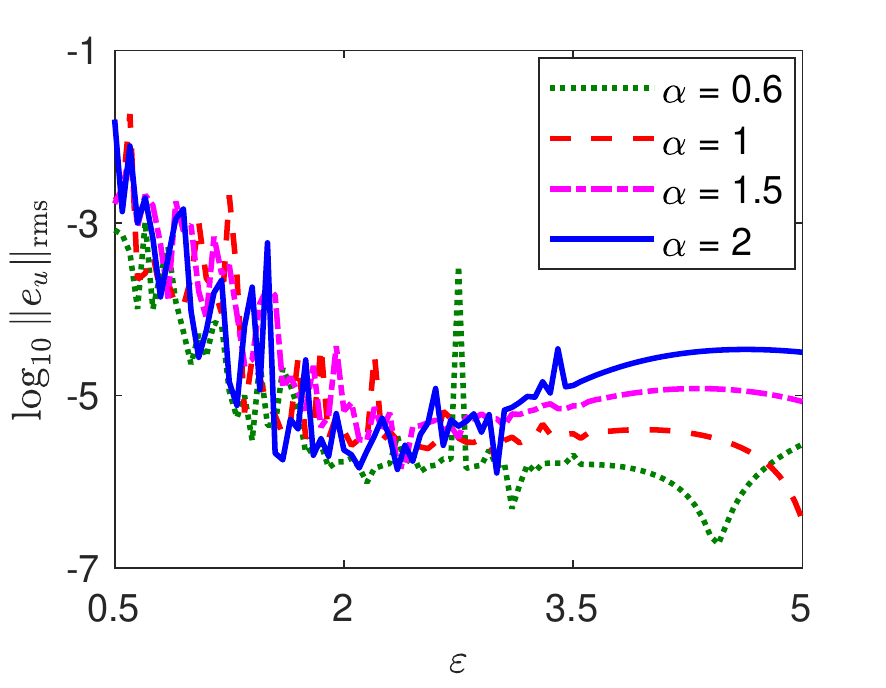}}
\centerline{(c)\includegraphics[height = 5.76cm, width = 7.86cm]{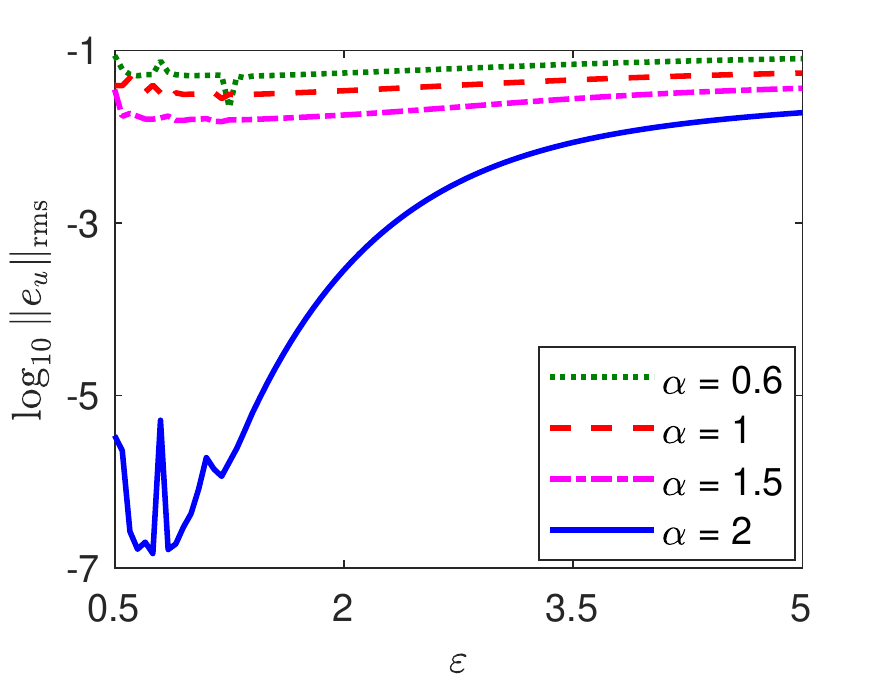}\hspace{-3mm}
(d)\includegraphics[height = 5.76cm, width = 7.86cm]{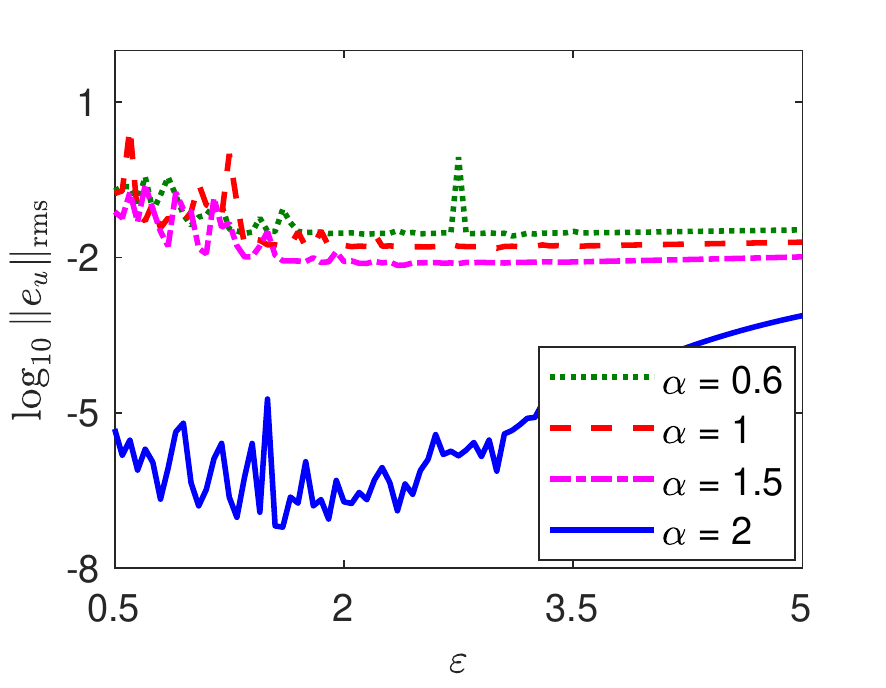}}
\caption{Numerical errors versus shape parameters in solving the 1D Poisson problem with $f$ and $g$ in (\ref{fun4-1})--(\ref{fun4-2}).  
(a) $\bar{N} = 17, s = 3$;  (b) $\bar{N} = 33, s = 3$;  (c) $\bar{N} = 17, s = 0$;   (d) $\bar{N} = 33, s = 0$.}\label{Fig4-1-2}
\end{figure}

As discussed previously,  the shape parameter $\varepsilon$ plays an important role in the accuracy of RBF-based methods.  
To study it, we compare the numerical errors of different shape parameters in Fig. \ref{Fig4-1-2}. 
It shows that the optimal shape parameter depends not only on exponent $\ap$ but also on solution regularity. 
For $s = 3$, large numerical errors are found if the shape parameter is  too small or too big. 
However, numerical errors for $s = 0$ and $\ap < 2$ become almost insensitive to the shape parameter  (see Fig. \ref{Fig4-1-2} (c) and (d)), suggesting that the solution regularity caps the numerical errors in this case. 
Moreover, the dependence of numerical errors on the shape parameter becomes more complicated when the number of RBF center points increases, and multiple optimal shape parameters might occur;  (cf. Fig. \ref{Fig4-1-2} (a) and (b)). 
We will leave the investigation of optimal shape parameters in solving nonlocal problems as our future research. 

Next, we compare our method with the recently proposed pseudospectral method in \cite{Rosenfeld2019},  where the benchmark fractional Poisson problem with $s = 0$ was  studied. 
The method in \cite{Rosenfeld2019} is different from ours mainly in the following aspects: 
\vspace{-1mm}
\begin{itemize}\itemsep -2pt
\item[(i)] It applies the extended Dirichlet boundary conditions to the pseudo-differential form of the fractional Laplacian in (\ref{pseudo}). 
To this end, this method requires a much larger computational domain than the physical domain $\Og$. 
For instance, a computational domain of $[-8, 8]$ was taken in solving the Poisson problem on $\Og = (-1, 1)$; see \cite[Section 4.2]{Rosenfeld2019}. 
This significantly increases the computational cost and storage requirement as more points are demanded in their computations. 
{ Fig. \ref{Fig4-1-3} (a) compares our numerical errors with those in \cite[Table 2]{Rosenfeld2019} in solving the 1D fractional Poisson problem with $s = 0$, where the uniform distance between center points is denoted as $h$.  
Numerical errors of these two methods are comparable, but the method  in \cite{Rosenfeld2019} used 8 times more points than ours due to its larger computational domain. 
}
\begin{figure}[htb!]
\centerline{(a)\includegraphics[height = 5.66cm, width = 7.86cm]{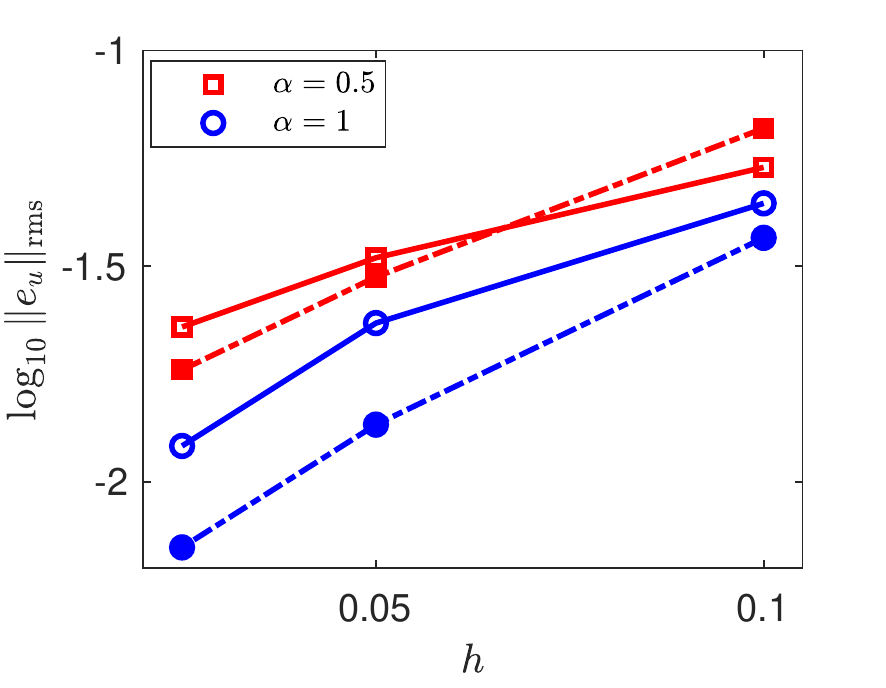}\hspace{-3mm}
(b)\includegraphics[height = 5.66cm, width = 7.86cm]{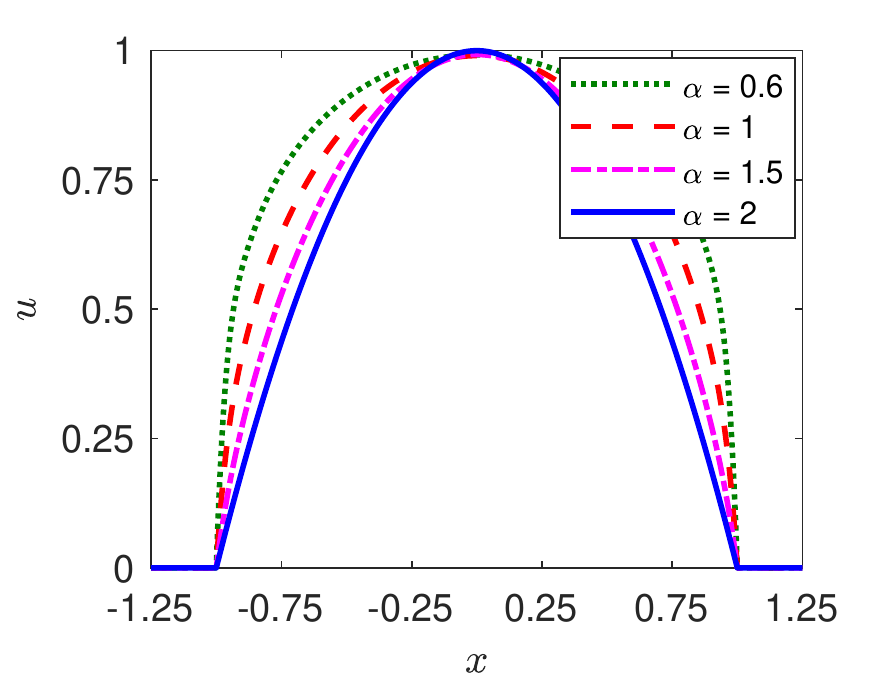}}
\caption{(a) Comparison of our numerical errors (dashed line) with those in \cite[Table 2]{Rosenfeld2019} (solid line) in solving 1D Poisson problem with $s = 0$. 
(b) Numerical solution of the 1D Poisson problem with $s = 0$, showing no Gibbs phenomenon of our method, where the shape parameter $\varepsilon = 4.3$.}\label{Fig4-1-3}
\end{figure}
\item[(ii)] Even though a larger computational domain is adopted, it has difficulties to handle the exact boundary conditions over $\Omega^c$. 
In fact, the boundary conditions only on a small region $\og \subset \Og^c$ is considered in their method; see more discussion in Remark \ref{remark1}. 
Consequently,  as pointed out in \cite[Figs. 2 and 3]{Rosenfeld2019}, the Gibbs phenomenon was observed near the boundary points $x = \pm 1$, and large numerical errors were found outside of the domain $[-1,1]$. 
In contrast, our method is free of these issues (see Fig. \ref{Fig4-1-3} (b)),  as it exactly utilizes the boundary conditions  and avoids the RBF approximation on ${\mathbb R}\backslash\bar{\Og}$. 
\item[(iii)] Different from ours, the method in \cite{Rosenfeld2019} uses the compactly supported Wendland RBFs as the basis function. Since the fractional Laplacian of  Wendland RBFs is unknown, numerical quadrature rules are required for their  approximation.  
This not only complicates the practical implementation of this method but hinders its generalization to problems  with the classical Laplacian. 
\end{itemize}

\subsection{Two-dimensional  Poisson problems}
\label{section4-2}

Our meshfree method can easily handle complex geometry, and its computational complexity is independent of  dimension $d \ge 1$. 
In this section, we will study its performance in solving the two-dimensional Poisson problem (\ref{BVP})--(\ref{BC}) on both regular and irregular domains. 
\bb

\subsubsection{Regular domain} 
\label{section4-2-1}

Here, we consider the two-dimensional  Poisson problem (\ref{BVP})--(\ref{BC}) on a unit disk domain, i.e., $\Og = \{\bx \in {\mathbb R}^2\,  |\, |\bx| < 1\}$. 
We choose function $f(\bx) \equiv 1$ for $\bx \in \Og$, and homogeneous Dirichlet boundary conditions $g(\bx) \equiv 0$ for $\bx \in\Upsilon$. 
In this case, the exact solution  is given by 
\beas
u(\bx) = 2^{-\ap}\Big(\Gamma(1+\fl{\ap}{2})\Big)^{-2}\big(1-|\bx|^2\big)^{\fl{\ap}{2}}, \qquad \mbox{for} \ \ \bx \in \Og,\quad \ap \in (0, 2].
\eeas
In our simulations, we choose the RBF center and test points radially distributed  on the disk domain $\bar{\Og}$. 
Choose an integer $n \ge 1$, and define the sets 
\beas
{\mathcal S}_{\bar{\Og}}^c = {\mathcal S}_{\bar{\Og}}^t = \bigg\{\fl{l}{n}\Big(\cos\big(2j\pi/(n+1)\big), \ \sin\big(2j\pi/(n+1)\big)\bigg), \ \ 
\mbox{for} \  0 \le l \le n, \ 0\le j \le n\Big\}. 
\eeas
that is, the total number of points is $\bar{N} = n(n+1)+1$ with $n$ radial layers. 

\begin{figure}[htb!]
\centerline{
\includegraphics[height = 3.86cm, width = 4.46cm]{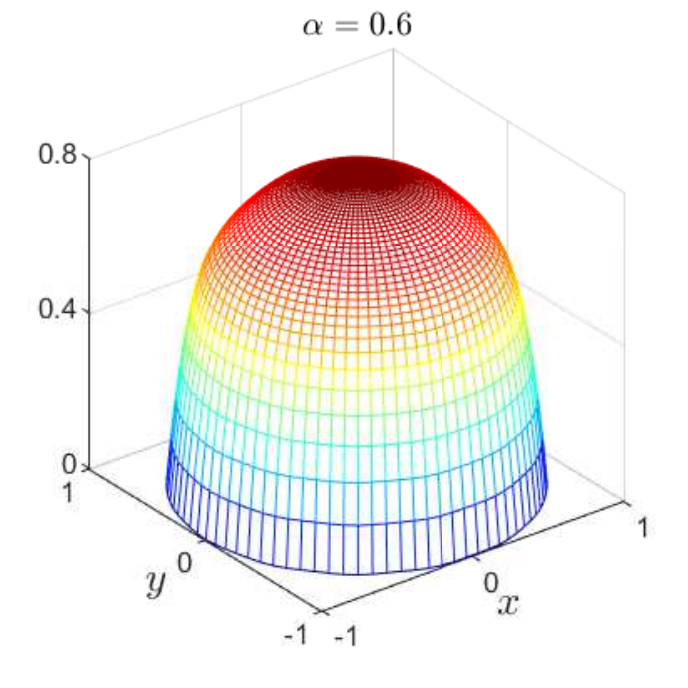}
\includegraphics[height = 3.86cm, width = 4.46cm]{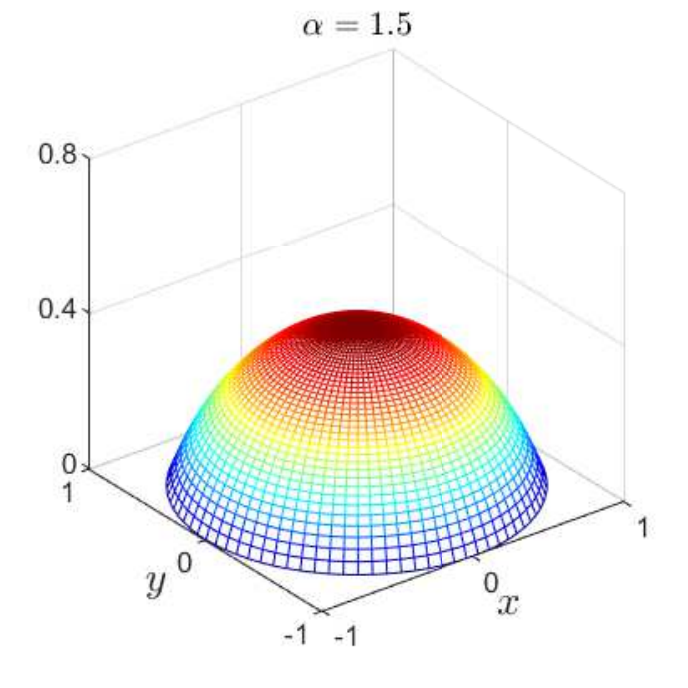}
\includegraphics[height = 3.86cm, width = 4.46cm]{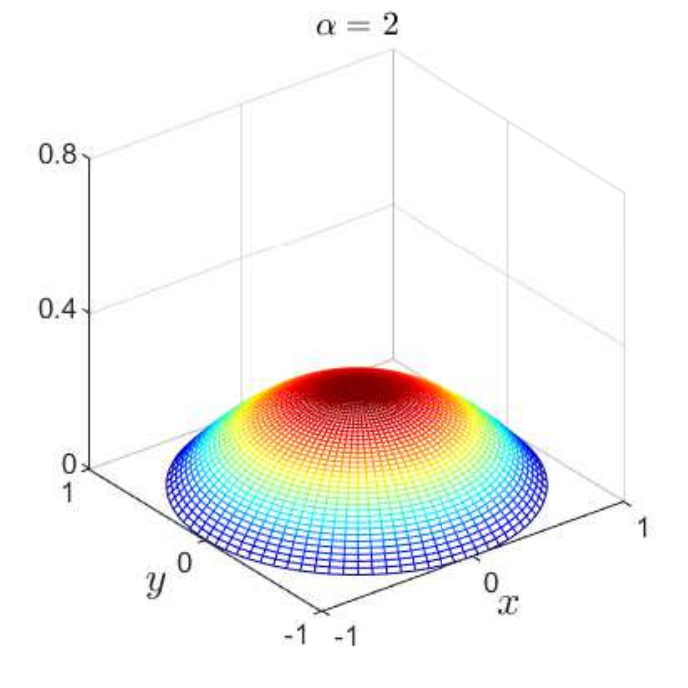}}
\caption{Numerical solution of the Poisson problem on a unit disk domain.}\label{Fig4-2-1}
\end{figure}

\begin{figure}[htb!]
\centerline{
(a)\includegraphics[height = 6.06cm, width = 6.36cm]{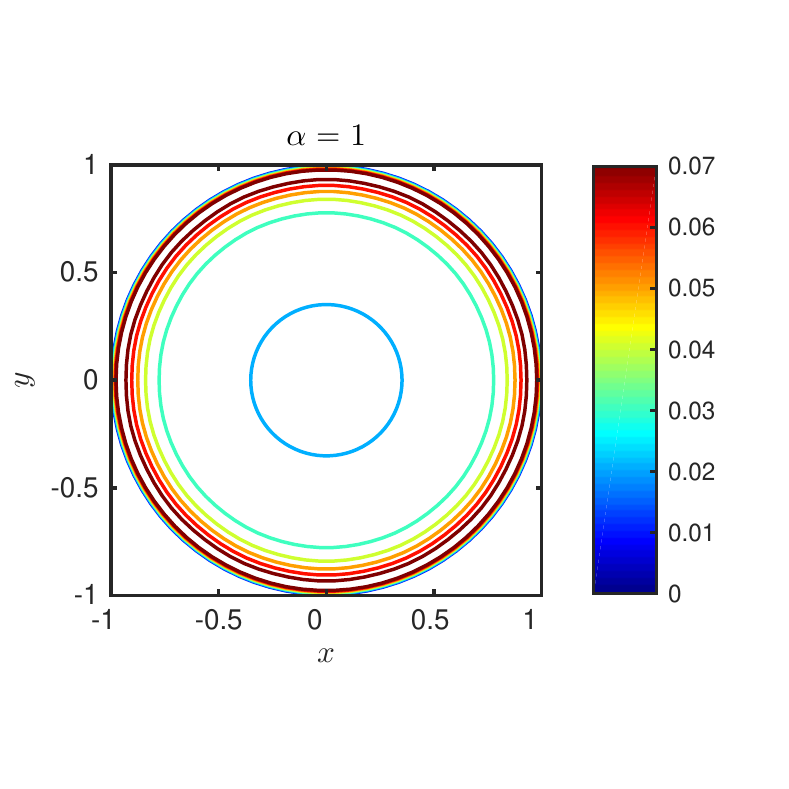}
(b)\includegraphics[height = 5.16cm, width = 7.56cm]{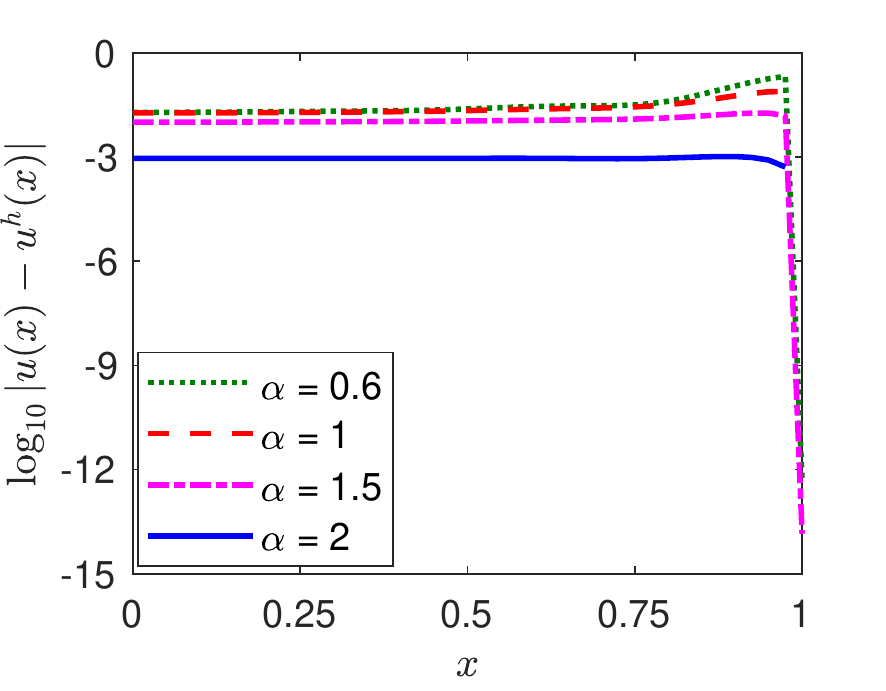}}
\caption{Pointwise numerical errors in solving the Poisson problem  on a unit disk domain, where $\ap = 1$ and $y = 0$ is fixed in (b). }\label{Fig4-2-2}
\end{figure}

In this example,  the solution has low regularity up to the boundary, i.e. $u \in C^{0, \fl{\ap}{2}}(\bar{\Og})$ \cite{Duo-FDM2019, Dyda2012, Acosta2017}. 
The smaller the exponent $\ap$, the less smooth the solution near boundaries; see Fig. \ref{Fig4-2-1} for numerical solution  of $\ap = 0.6, 1.5$  and $2$. 
Moreover, the solution becomes much flatter as $\ap$ increases. 
Table \ref{Tab4-2-1} presents numerical errors $\|e_u\|_{\rm rms}$ for different  $n$, where we take the shape parameter  $\veps = 2$. 
It shows that { as the number of points increases}, the numerical errors  { decrease with a spectral rate for $\ap = 2$.  
While $\ap < 2$,  the solution regularity dominates the problem, and numerical errors decrease slowly, which is similar to the observations in Table \ref{Tab4-1-1} for $s = 0$. 
}
\begin{table}[htb!]
\begin{center} 
\begin{tabular}{|c||c|c|c|c|c|}
\hline
 $n$ & $\ap = 0.6$ & $\ap = 1$ & $\ap = 1.5$ & $\ap  = 2$\\
\hline
$3$  &1.360E-1& 8.187E-2& 3.817E-2& 1.211E-2\\
\hline
$4$ &9.802E-2& 5.640E-2& 2.459E-2& 5.924E-3\\
\hline
$5$ &7.471E-2& 4.130E-2& 1.667E-2& 2.488E-3\\
\hline
$6$ &5.925E-2& 3.178E-2& 1.204E-2& 9.059E-4\\
\hline
$7$ &4.844E-2& 2.544E-2& 9.229E-3& 2.803E-4\\
\hline
\end{tabular}
\caption{Numerical errors $\|e_u\|_{\rm rms}$ in solving the Poisson problem on a unit disk domain, where the shape parameter  $\varepsilon = 2$.}\label{Tab4-2-1}
\end{center}
\end{table}
Fig. \ref{Fig4-2-2} further shows the pointwise errors for $\ap = 1$. 
We find that numerical errors are radially symmetric and the maximum error is found around the domain boundary.  
This suggests that the local refinement  (including more points around boundary) might improve the accuracy of this problem.  
Our extensive studies show that numerical errors are considerably smaller if the solution is smoother around boundary. 

Compared to other methods, the geometric flexibility enables our method to solve this problem with much less computational efforts.  
For instance of $\ap = 1$, our method can achieve the accuracy of ${\mathcal O}(10^{-2})$ with  the total number of points $\bar{N} \approx 10$, but the Wendland RBF method requires at least 400 points { due to its larger computational domain} (see \cite[Table 4]{Rosenfeld2019}).  
On the other hand, finite element methods (FEM) are also flexible of  domain geometry, but they are mesh-based local approximation methods. 
The results in \cite{Bonito2019} show that to achieve the same accuracy, much larger number of unknowns are introduced by FEMs; see \cite[Table1]{Bonito2019}. 
Moreover, the nonlocality and singularity of the fractional Laplacian makes assembling the stiffness matrix  extremely challenging for FEMs \cite{Acosta2017-2D, Bonito2019}.

\subsubsection{Irregular domain}
\label{section4-2-2}

In this section, we solve the two-dimensional Poisson problem  (\ref{BVP})--(\ref{BC})   on an irregular domain $\Og = \big\{\bx\,|\, \bx \in (-1, 1)^2\ {\rm and} \ |\bx| >  0.5\big\}$, i.e. a two-dimensional domain confined between a unit square and a circle with radius $r = 0.5$.
We will consider an inhomogeneous Dirichlet boundary condition and choose
\beas
\label{2D-Poisson}
f(\bx) = \Gamma(2+\ap) \,_2F_1\Big(\fl{2+\ap}{2}, \fl{3+\ap}{2}; 1; -|\bx|^2\Big),  && \quad\mbox{for} \  \ \bx \in \Og, \qquad\\
\label{2D-BC}
g(\bx) = \sqrt{(1+|\bx|^2)^{-3}}, &&\quad\mbox{for} \ \ \bx \in  \Upsilon.
\eeas
In this case, the exact solution of the Poisson problem can be constructed as $u(\bx) = \sqrt{(1+|\bx|^2)^{-3}}$ for $\bx \in {\mathbb R}^2$ and $\ap \in (0, 2]$. 

In our simulations, we take the shape parameter $\varepsilon = 1.5$. 
The test points are chosen to be the same as center points.
To this end, we first choose a set of points radially distributed on an annulus for $0.5 \le |\bx| \le 1$, i.e., for $n \in {\mathbb N}$
\beas
\hat{\mathcal S} = \bigg\{\Big(\fl{1}{2}+\fl{l}{2n}\Big)\Big(\cos\big(j\pi/2n\big), \ \sin\big(j\pi/2n\big)\bigg), \ \ 
\mbox{for} \  0 \le l \le n, \ 1\le j \le 4n\Big\}. 
\eeas
Then we map the points in $\hat{\mathcal S}$ to domain $\Og$ using the elliptic grid mapping \cite{Fong0014},  i.e., letting 
\beas
&&x_i = \fl{1}{2} \bigg(\sqrt{2+\hat{x}_i^2-\hat{y}_i^2+2\sqrt{2}\hat{x}_i} - \sqrt{2+\hat{x}_i^2-\hat{y}_i^2-2\sqrt{2}\hat{x}_i}\bigg), \\
&&y_i = \fl{1}{2} \bigg(\sqrt{2-\hat{x}_i^2+\hat{y}_i^2+2\sqrt{2}\hat{y}_i} - \sqrt{2-\hat{x}_i^2+\hat{y}_i^2-2\sqrt{2}\hat{y}_i}\bigg), 
\eeas
for $(\hat{x}_i, \hat{y}_i) \in \hat{\mathcal S}$.
That is,  the total number of points is $\bar{N} = 4n(n+1)$.

\begin{figure}[htb!]
\centerline{(a) \includegraphics[height = 4.86cm, width = 5.26cm]{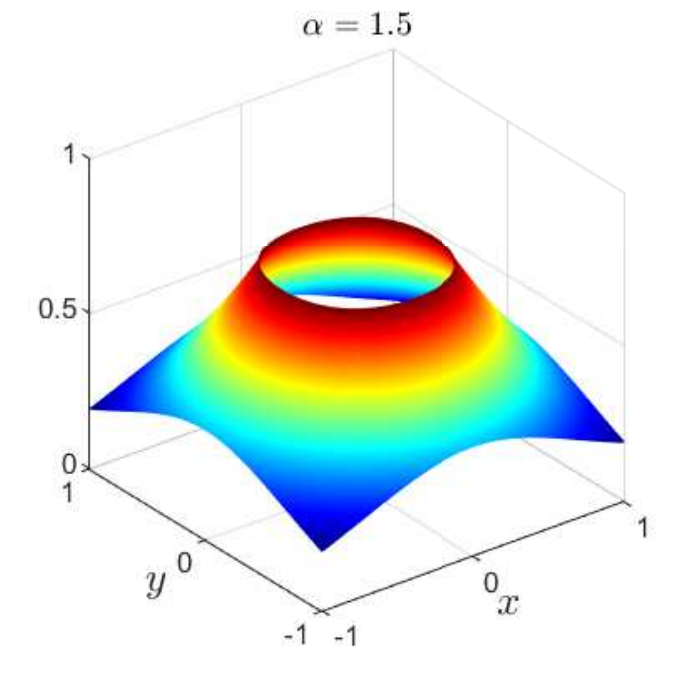}
(b) \includegraphics[height = 4.86cm, width = 6.06cm]{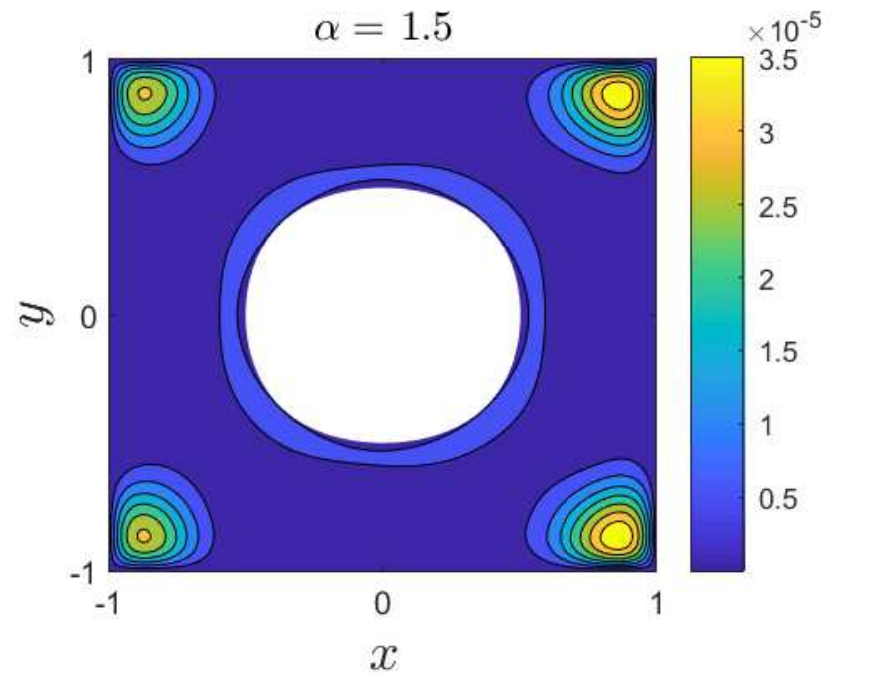}}
\caption{Numerical solution (a) and numerical errors (b) in solving the two-dimensional Poisson problem on an irregular domain, where we choose $n = 5$ in our simulations.}\label{Fig4-2-3}
\end{figure}

Fig. \ref{Fig4-2-3} shows the numerical solution and pointwise errors for $\ap = 1.5$. 
We find that numerical errors reach the maximum at four corners of the domain, but overall they are  small even with the number of points $\bar{N} = 24$. 
Table \ref{Tab4-3} further demonstrates the numerical errors $\|e_u\|_{\rm rms}$ { and condition number ${\mathcal K}$ of the linear system} for various $\ap$. 
\begin{table}[htb!]
\begin{center}
\begin{tabular}{|c||c|c|c|c|c|c|c|c|}
\hline
& \multicolumn{2}{|c|}{$\ap = 0.6$} & \multicolumn{2}{|c|}{$\ap = 1$} & \multicolumn{2}{|c|}{$\ap = 1.5$} & \multicolumn{2}{|c|}{$\ap = 2$}   \\
\cline{2-9}
 n & $\|e_u\|_{\rm rms}$ & $\mathcal{K}$ & $\|e_u\|_{\rm rms}$ & $\mathcal{K}$ & $\|e_u\|_{\rm rms}$ & $\mathcal{K}$ & $\|e_u\|_{\rm rms}$ & $\mathcal{K}$  \\ 
\hline
$2$  &1.635E-3 & 2.313E7  & 1.772E-3 & 3.143E7  &2.044E-3 & 4.694E7 &  2.503E-3 &6.955E7\\
\hline
$3$  &2.983E-4 & 3.50E11 & 3.958E-4 & 5.69E11 &6.156E-4 & 1.03E12 & 1.006E-3 &1.51E12\\
\hline
$4$  &8.118E-5 & 1.69E16 & 9.821E-5 & 2.48E16 &1.554E-4 & 4.58E16 & 9.821E-5 &4.41E16\\
\hline
\end{tabular}
\caption{Numerical errors $\|e_u\|_{\rm rms}$ and condition number ${\mathcal K}$ in solving the two-dimensional Poisson problem on an irregular domain, where the shape parameter $\varepsilon = 1.5$.}\label{Tab4-3}
\end{center}
\end{table}
Compared to the results in Table \ref{Tab4-2-1},  numerical errors reduce much faster as the number of points increases, since the solution in this case are infinitely differentiable on ${\mathbb R}^2$. 
{ With $\bar{N}$ further increasing, the system may become ill-conditioned and affect the simulation errors. 
Recently,  multi-precision toolboxes and domain decompositions have been proposed to resolve this issue; see \cite{Kansa2017, Sarra2011, Sarra2017} and references therein for more discussion.
}

\subsection{Diffusion problems}
\label{section4-3}

In this section, we further study the performance of our method in solving time-dependent problems. 
To this end, we consider the following diffusion problem with nonhomogeneous Dirichlet boundary conditions: 
\bea\label{diffusion}
\p_t u(\bx, t) = -(-\Dt)^\fl{\ap}{2}u + f(\bx, t),  & \ \ & \mbox{for} \ \ \bx\in \Og =  (-1, 1)^2, \quad t > 0, \qquad  \\
\label{BC-new}
u(\bx, t)  = t(1+0.5|\bx|^2)^{-\fl{3}{2}},  && \mbox{for} \ \ \bx\in \Upsilon, \quad t \ge 0, \\
\label{initial}
u(\bx, 0) = 0, && \mbox{for} \ \ \bx\in {\mathbb R}^2.
\eea
The source term $f$ is chosen such that the exact solution of (\ref{diffusion})--(\ref{initial}) is given by:
\bea\label{exact4-3}
u(x, y) = t /\sqrt{\big(1 + 0.5(x^2 + y^2)\big)^{3}}, \qquad \mbox{for} \ \ \bx\in {\mathbb R}^2, \quad t \ge 0. 
\eea

Denote time sequence $t_n = n\Dt t$ (for $n = 0, 1, \ldots$) with time step $\Dt t > 0$. 
Choose the RBF center and test points as uniformly distributed tensor grid points on $[-1, 1]^2$. 
We then discretize  (\ref{diffusion}) with  the Crank--Nicolson method in time and our meshfree method in space and obtain the fully discretized scheme as:
\bea\label{scheme}
\fl{u^{n+1}_k - u^{n}_k}{\Dt t} = -(-\Dt)_h^\fl{\ap}{2} \Big(\fl{u_k^n + u_k^{n+1}}{2}\Big) + \fl{f(\bx_k, t_n) + f(\bx_k, t_{n+1})}{2}, \quad \mbox{for} \ \ 
\bx_k \in {\mathcal S}_\Og^t, \ \ n = 0, 1, \ldots,
\eea 
where $(-\Dt)_h^\fl{\ap}{2}u_k^n$ represents the numerical approximation of $(-\Dt)^\fl{\ap}{2}u(\bx_k, t_n)$ with $u_k^n$ being the approximation of $u(\bx_k, t_n)$, i.e.\bea\label{anzta}
u^n_k= \sum_{1 \le i \le \bar{N}} \lambda_i^n \varphi^\veps(|\bx_k - \bx_i|). 
\eea
At each time step, we will solve for the unknowns $\lambda_i^{n+1}$ for $1 \le i \le \bar{N}$. 
Simplifying (\ref{scheme}) and applying (\ref{anzta}) to boundary points $\bx_k \in {\mathcal S}_{\p\Og}^t$, we then obtain our scheme in matrix-vector form: 
\bea\label{scheme0}
\left(\begin{array}{c}\displaystyle B_{N\times\bar{N}} + \fl{\Dt t}{2} A_{N\times\bar{N}}\\
\displaystyle C_{(\bar{N}-N)\times\bar{N}}\end{array}\right)
\Lambda_{\bar{N}\times 1}^{n+1} = 
\left(\begin{array}{c}\displaystyle \Big(B - \fl{\Dt t}{2} A\Big)\Lambda^n+ \fl{\Dt t}{2}\big({\bf f}^n + {\bf f}^{n+1} + {\bf w}^n + {\bf w}^{n+1}\big) \\
\displaystyle{\bf g}^{n+1}\end{array}\right), 
\eea
for $n = 0, 1, \ldots$, where we denote $\Lambda_{\bar{N}\times 1}^n = \big(\lambda_1^n,  \cdots, \lambda_{\bar{N}}^n\big)^T$,  \ ${\bf f}_{N\times1}^n = \big(f(\bx_1, t_n),  \cdots,  f(\bx_N, t_n)\big)^T$,   and ${\bf g}_{(\bar{N}-N)\times 1}^n = \big(g(\bx_{N+1}, t_n), \cdots, g(\bx_{\bar{N}}, t_n)\big)^T$ with $g$ denoting the boundary condition in (\ref{BC-new}). 
The vector ${\bf w}_{N\times1}^n = \big(w(\bx_1, t_n),  \cdots,  w(\bx_N, t_n)\big)^T$ is from the discretization of $(-\Dt)^\fl{\ap}{2}$ in (\ref{Eq-uniform}) with
\beas
w(\bx, t) = \zeta_\ap C_{d, \ap} \int_{\Og^c} \fl{g({\bf y}, t)}{|\bx - \by|^{d+\ap}} d{\bf y}, 
\eeas
which reduces to zero if $\ap = 2$ or homogeneous boundary conditions are considered. 
The matrices $B$ and $C$ are composed of  the Gaussian RBFs $\varphi^\veps(|\bx_k - \bx_i|)$, while the entries of matrix $A$ are given by the coefficients of $\lambda_i$ in (\ref{Eq-gov}).
All of these matrices  remain the same at all time steps.

\begin{figure}[htb!]
\centerline{
\includegraphics[height = 3.86cm, width = 4.46cm]{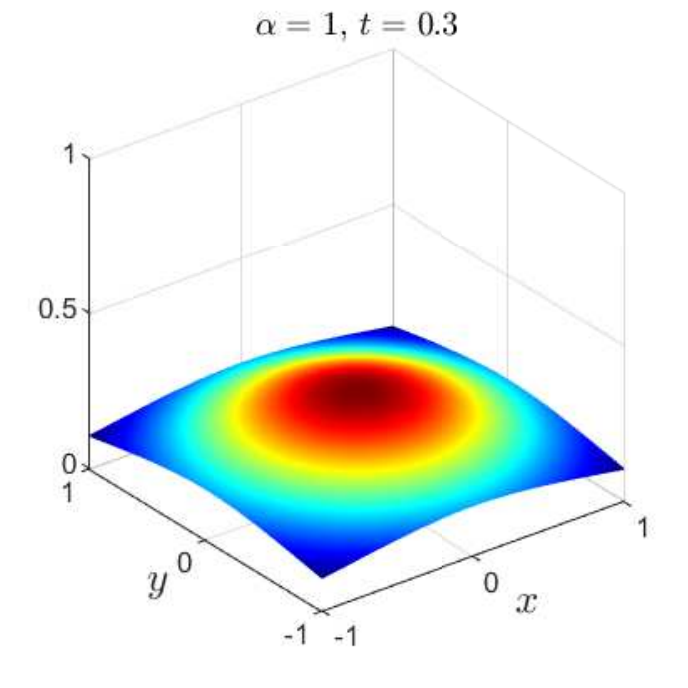}
\includegraphics[height = 3.86cm, width = 4.46cm]{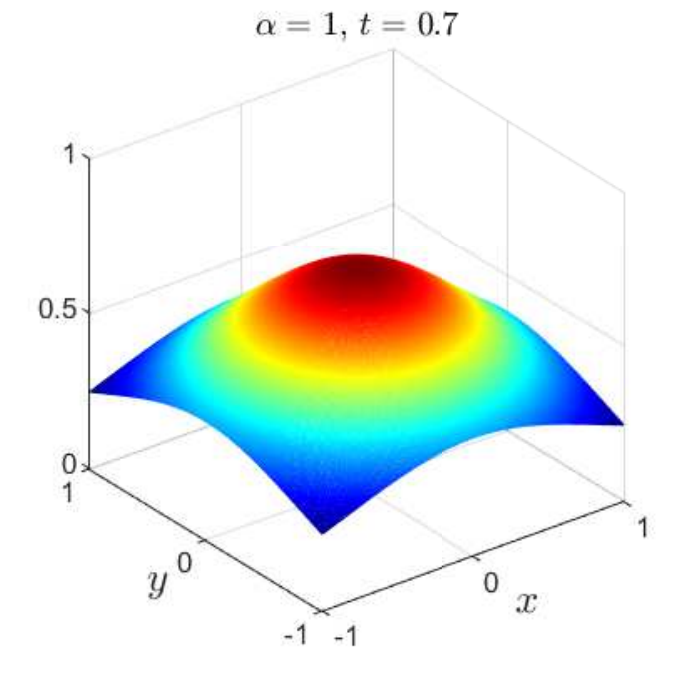}
\includegraphics[height = 3.86cm, width = 4.46cm]{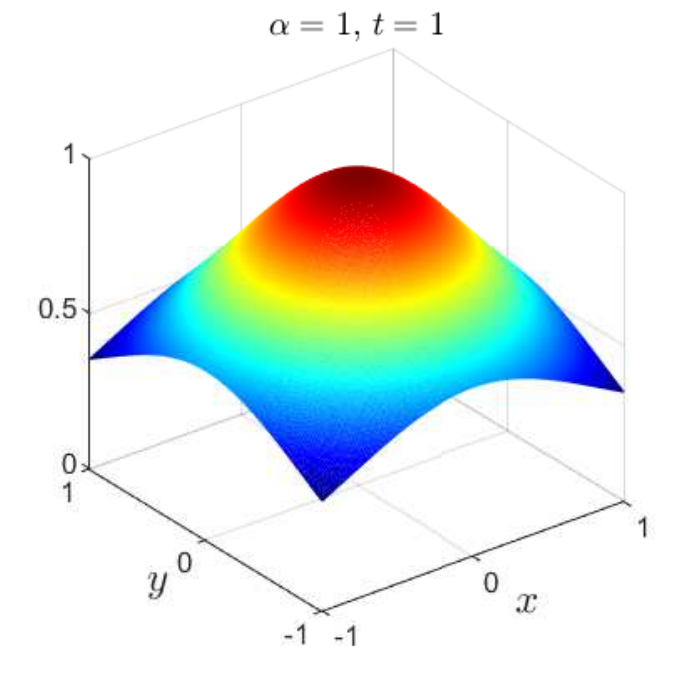}}
\caption{Time evolution of the solution of the diffusion problem (\ref{diffusion})--(\ref{initial}).}\label{Fig4-3-1}
\end{figure}

\begin{figure}[htb!]
\centerline{
(a)\includegraphics[height = 5.76cm, width = 7.86cm]{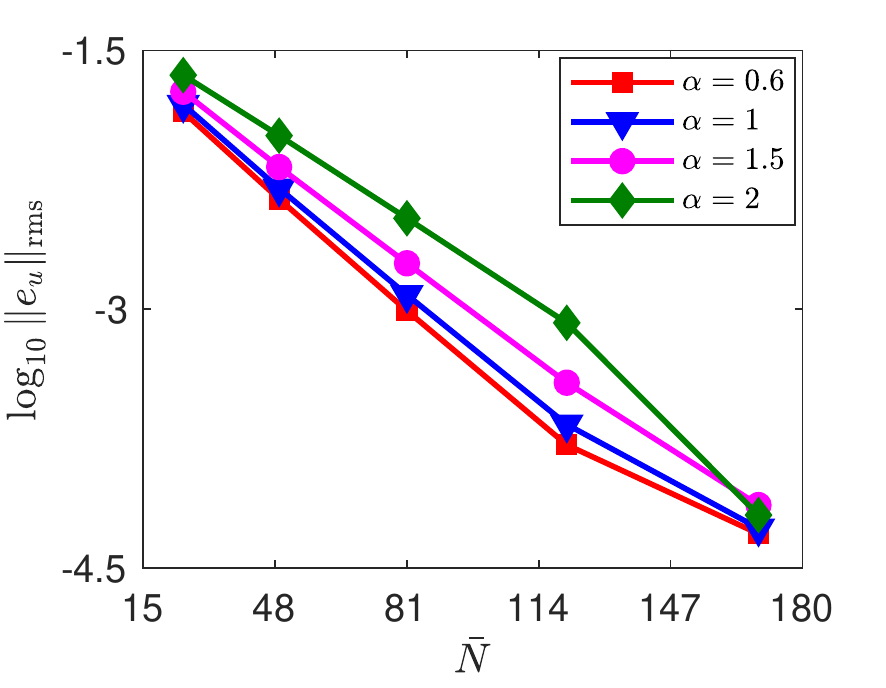}
\hspace{-5mm}
(b)\includegraphics[height = 5.76cm, width = 7.86cm]{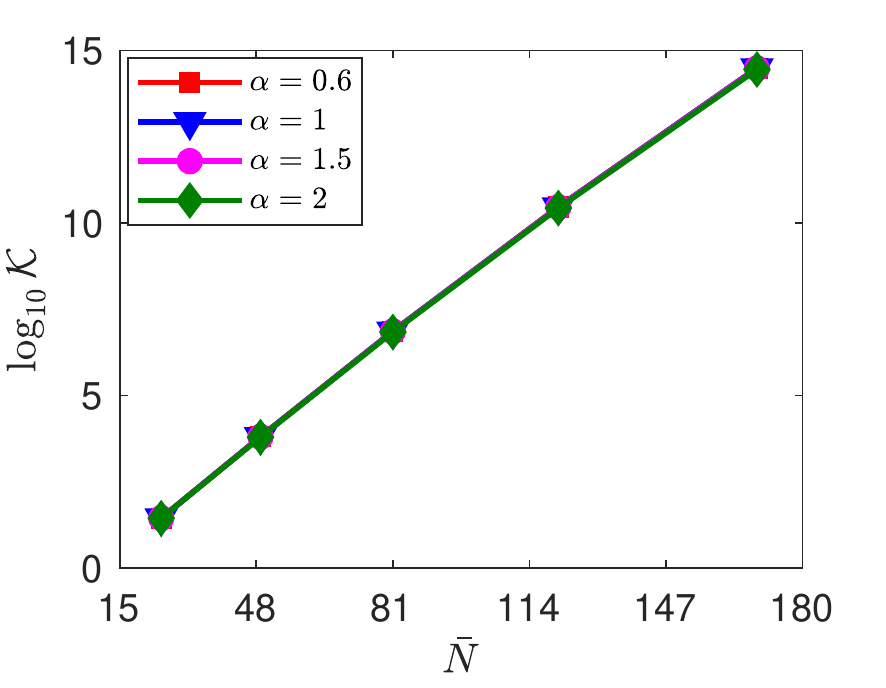}}
\caption{Numerical errors $\|e_u\|_{\rm rms}$ at time $t = 1$ and condition number ${\mathcal K}$ in solving the diffusion problem (\ref{diffusion})--(\ref{initial}), where the shape parameter $\veps = 1.9$.}\label{Fig4-3-2}
\end{figure}

In our simulations, we take a small time step $\Dt t=0.001$ such that the temporal errors are neglectable in comparison  to spatial errors. 
Fig. \ref{Fig4-3-1} shows the time evolution of the solution for $\ap = 1$, which agree well with the exact solution in (\ref{exact4-3}). 
Furthermore, we present the numerical errors $\|e_u\|_{\rm rms}$ at $t = 1$ and the condition number ${\mathcal K}$ of the linear system in Fig. \ref{Fig4-3-2}. 
It shows that numerical errors decrease with a spectral rate for any $\ap\in(0, 2]$, while the condition number increases quickly as more points are used. 
If the number of points $\bar{N}$ is too big, one can adopt the strategies in \cite{Sarra2011, Kansa2017, Sarra2017} to avoid the ill-conditioning issue. 

\section{Summary and discussions}
\setcounter{equation}{0}
\label{section5}

We proposed a novel meshfree pseudospectral method to solve both the classical and fractional PDEs in a unified scheme. 
Our method takes great advantage of the Laplacian of the Gaussian RBFs and enables us to approximate the classical and fractional Laplacian in a single framework, which is the key merit of our method distinguishing it from other existing methods for fractional PDEs. 
Extensive numerical experiments were carried out to study the performance of our method in  approximating the Dirichlet Laplace operators and solving the classical and fractional PDE problems. 
Compared to mesh-based methods, our method can easily handle complex geometry and achieve higher accuracy with fewer points. 
More importantly, our method can  solve the $d$-dimensional (for $d \ge 1$) classical and fractional PDEs with a single computer implementation, which could greatly benefit the study of coexistence of normal and anomalous diffusion in recent applications. 
We compared our method with the recently proposed Wendland RBF method in \cite{Rosenfeld2019}. 
In contrast to it, our method  exactly incorporates the Dirichlet boundary conditions into the scheme and is free of the Gibbs phenomenon as observed in \cite{Rosenfeld2019}.  
Moreover, the method in \cite{Rosenfeld2019} solves the problem on a computational domain that is much larger than the physical domain $\bar{\Og}$,  and consequently its computational complexity is much larger than ours. 
Our studies suggested that to obtain good accuracy the shape parameter cannot be too small or too big, and the optimal shape parameter depends on the RBF center points,  the exponent $\ap$, and the solution properties. 
How to find the optimal shape parameter is still an active research topic  in the area of RBF-based methods. 
We will leave it for our future study, especially in solving nonlocal problems.

\bigskip
\noindent{\bf Acknowledgements. } 
We acknowledge two anonymous reviewers for their valuable comments that greatly help to improve  this manuscript. 
This work was supported by the US National Science Foundation under Grant Number DMS-1913293 and DMS-1953177. 

\end{document}